\documentclass[12pt]{article}
\usepackage{color}
\usepackage{amsfonts}
\usepackage{amssymb}
\usepackage{amscd}
\usepackage{amsmath}
\usepackage{epsfig}
\usepackage{graphicx, subfigure}
\usepackage{latexsym}
\usepackage{pst-3dplot}
\usepackage{verbatim}
\usepackage{color}

\usepackage{color, graphics}
\usepackage{float}
\usepackage{tikz}
\usepackage{pdfsync}
\usepackage{hyperref}
\usepackage{multirow}

\usepackage{calrsfs}

\newlength{\dinwidth}
\newlength{\dinmargin}
\newtheorem{theorem}{Theorem}

\def\diag{\mathrm{diag}}

\def\px1{p_{x_1}}
\def\px2{p_{x_2}}
\def\pu1{p_{u_1}}

\begin{document}
\title{Demchenko's nonholonomic case of gyroscopic ball rolling without sliding over a sphere after his 1923 Belgrade doctoral thesis}

\author{Vladimir Dragovi\'c$^1$, Borislav Gaji\'c$^2$ and Bo\v zidar Jovanovi\'c$^3$}

\date{}

\maketitle

\footnotetext[1]{Department of Mathematical Sciences, University
	of Texas at Dallas, 800 West Campbell Road, Richardson TX 75080,
	USA. Mathematical Institute SANU, Kneza Mihaila 36, 11000
	Belgrade, Serbia.  E-mail: {\tt
		Vladimir.Dragovic@utdallas.edu}--the corresponding author}

\footnotetext[2]{Mathematical Institute SANU, Kneza Mihaila 36, 11000
	Belgrade, Serbia.  E-mail: {\tt
		gajab@mi.sanu.ac.rs}}

\footnotetext[3]{Mathematical Institute SANU, Kneza Mihaila 36, 11000
	Belgrade, Serbia.  E-mail: {\tt
		bozaj@mi.sanu.ac.rs}}

\begin{abstract}
	We present an integrable nonholonomic case of rolling without sliding of a gyroscopic ball over a sphere.  This case was introduced and studied in detail by Vasilije Demchenko in his 1923 doctoral dissertation defended at the University of Belgrade, with Anton Bilimovi\'c as the advisor. These results are absolutely unknown to modern researchers. The study is based on the C. Neumann coordinates and the Voronec principle. By using involved technique of elliptic functions, a detailed study of motion is performed.
Several special classes of trajectories are distinguished, including regular and pseudo-regular precessions. So-called remarkable trajectories, introduced by Paul Painlev\'e and Anton Bilimovi\' c, are described in the present case. The historic context as well as the place of the results in contemporary mechanics are outlined.

	\vskip 1cm
	
	
	
\end{abstract}

\

{\it Dedicated to Anton Bilimovi\' c and his scientific school on the occasions of the 50 anniversaries of deaths
of Anton Bilimovi\' c and Vasilije Demchenko.}

\newpage

\maketitle
\tableofcontents

\section{Introduction}

We present an integrable nonholonomic case of rolling without sliding of a gyroscopic ball over a sphere.  This case was introduced and studied in detail by Vasilije Demchenko in his 1923 doctoral dissertation defended at the University of Belgrade, with Anton Bilimovi\'c  as the advisor. These results are absolutely unknown to modern researchers. The study is based on the C. Neumann coordinates and the Voronec principle. By using involved technique of elliptic functions, a detailed study of motion is performed.
Several special classes of trajectories are distinguished, including regular and pseudo-regular precessions. So-called remarkable trajectories, introduced by Paul Painlev\'e and Anton Bilimovi\' c, are described in the present case. The historic context as well as the place of the results in contemporary science are outlined.

\medskip

Anton Bilimovi\'c (1879-1970) was an outstanding representative of the Russian scientific elite, who ended up in Belgrade
as a result of the turmoil induced by the revolution in Russia. As an already established scientist and a former Rector of the Malorosiisk University in Odessa, Bilimovi\'c made a tremendous contribution to the further development and the organization of mathematics and mechanics in Belgrade, Serbia, and Yugoslavia. A detailed biography of Academician Bilimovi\' c can be found in \cite{DjDj}. A comprehensive study of the scientific school of Anton Bilimovi\' c in Serbia till mid 1970's was given in \cite{GrFra}. It is presented there as an integral part of G. K. Suslov school. Namely,
 Bilimovi\' c's advisor was Peter Vasilievich Voronec (1871--1923), a distinguished pupil of Gavril Konstantinovich Suslov (1857--1935). \cite{GrFra} also describes the results of Serbian pupils of the Bilimovi\' c school: Tatomir Andjeli\'c (1903-1993), Rastko Stojanovi\'c (1926-1972), Veljko Vuji\v ci\'c (1929), Bo\v zidar Vujanovi\'c (1930-2014), and Djordje Djuki\'c (1943-2019). In his editorial concluding remarks in \cite{GrFra}, Andjeli\'c also listed works of
 Vencilsav Zhardecki (1896-1962) and Djordje Mu\v sicki (1921-2018) as parts of the school.

\medskip

Bilimovi\'c was not an isolated example of a Russian scientist who came to Belgrade at that time. Along with him, a notable scientist and his former scientific advisor, P. V. Voronec also came to Belgrade. According to one of the most romantic Belgrade urban stories, Peter and his son, Konstantin Voronec (1902-1974), met for the last time as recruited solders of different unites of "Whites".  Feeling that the civil war is not developing in the favor of "Whites" and that the end of the war is coming closer, they agreed upon to meet in Belgrade, with the condition that who comes first would wait for the other one every day at noon in front of the National Theatre. Peter, the father, reached Belgrade the first. Waiting for his sone day after day, he lost patience and decided to go back  and search for Konstantin. He requested Bilimovi\'c to take over his duty, and await for Konstantin in Belgrade. Peter and Konstantin did not manage to find each other. Konstantin on foot went across Romania, and came to Belgrade to be met and taken care of by Bilimovi\'c. At the same time, Peter ended up in Ukraine, got ill and died in 1923. The communications between two countries were so poor at that time, that Konstantin learned about his father's fate only after he came to Paris in 1930's. A comprehensive description of this dramatic period of the Voronec family was presented in \cite{SDj}.

\medskip

Along with the son, Bilimovi\'c, also took care about a pupil of Peter Voronec, Vasilije Gregorevich Demchenko (1898-1972). In 1923, Demchenko prepared and defended the doctoral dissertation:

\medskip

\emph{Rolling without sliding of a gyroscopic ball over a sphere}, University of Belgrade, 1923 (in Serbian) pp. 94.

\medskip

The dissertation was published as a separate book in Belgrade in 1924, see \cite{Dem1924}. The work was  motivated and based on the results of his teachers, Peter Voronec and Anton Bilimovi\'c, see \cite{Vor1902, Vor1912} and \cite{Bilim}. The dissertation was accepted for the final doctoral examination
on the meeting of the Faculty of Philosophy of the University of Belgrade, of 15 November of 1923, based on the report of the members of the examination committee: Anton Bilimovi\'c, Mihailo Petrovi\'c (1868-1943), and Milutin Milankovi\'c (1879-1952). Petrovi\'c was a founding father of the modern Serbian mathematics, see for example \cite{DG1, DG2, Petrovic}. Milankov\'c is best known for developing his mathematical theories of climate which, in particular produced one of the most significant theories relating Earth motions and long-term climate change, nowadays called Milankovi\'c cycles \cite{Mil1941} and \cite{Mil1969}.

\medskip

The dissertation consists of a Preface, six chapters and a Resume in French. In the last sentence of the Preface, Demchenko expresses the deepest gratitude to his teachers, Professors P. Voronec and A. Bilimovi\'c.

\medskip

The material related to nonholonomic problems of rolling of a surface over another surface accumulated
by that moment, Demchenko classifies in three categories: rolling of a ball over a surface, rolling of a surface over a plane, and rolling of a surface over the sphere. The case studied in the dissertation relates to the last category. It is interesting since it reduces to elliptic quadratures, as was indicated in\cite{Vor1912}. In a sense, this work is also continuation of works of Bobilev and Zhukovsky, \cite{Bob1892} and \cite{Zhuk1893}, who considered special cases of rolling of a gyroscopic ball over the plane.

\section{Demchenko's results: gyroscopic ball rolling without sliding over a sphere}
\subsection{The dissertation chapter by chapter}

The subject of the dissertation is at the interface of differential geometry of curves and surfaces, nonholonomic mechanics and the theory of elliptic functions and integrals. It is written in a clear and illuminating fashion, demonstrating author's full expertise in each of these fields and mastery in their synergy. The study is detailed and very well rounded. The obtained results are complete, numerous, interesting, transparent and rigorous. The exposition is elegant, with a well-thought organization which connects various chapters and subchapters into a fully focused and convergent material. There is a perfect measure between each detail and the global line as a whole. Let us also observe that the dissertation is written in a clean and smooth Serbian language,
with a few instances of constructions which could be seen as more natural in Russian than in Serbian.

\medskip

Chapter 1, Kinematics of rigid body rolling over a fixed surface, consists of three subchapters,
Section 1.1 Motion over the surface of Darboux trihedral; Section 1.2 Kinematic elements of a rolling rigid body, in Neumann coordinates; Section 1.3 The case of rolling without sliding.

\medskip

Chapter 2, The equations of motion of a rigid body, in a  moving frame with an arbitrary motion
with respect to the rigid body, has four subchapters. Section 2.1 The equations of motion of  a free rigid body in a moving frame; Section 2.2 The equations of motion of  a non-free rigid body; Sections 2.3 Applications to rolling without sliding; Section 2.4 Particular cases.

\medskip

Chapter 3, Voronec Principle, has three subchapters: 3.1 A principle similar to Hamiltonian, which is applicable to nonholonomic systems; Section 3.2 Application to rolling without sliding over a fixed surface; Section 3.3 Rolling of gyroscopic bodies.

\medskip

Chapter 4, Reducing to quadratures, consists of the following subchapters: Section 4.1 Bobilev problem and its generalization; Section 4.2 Kinematic elements and expressions for the kinetic energy; 4.3 Differential equations of motion and first integrals; Section 4.4 Calculation of coordinates $u$ and $v$; Section 4.5 Calculation of cyclic coordinates $u_1$, $v_1$, and $\vartheta$; Section 4.6 A particular Solution.

\medskip

Chapter 5, Solution in the finite form, has the following subchapters: Section 5.1 Inversion of the elliptic integral. Discriminant; Section 5.2 Arguments $a_0, b_0, a, b$; Section 5.21 Calculation of $v$;
Section 5.22 Calculation of $s, \tau$; Section 5.3 Arguments $a_1$ and $b_1$; Section 5.31 Calculation of $v_1$; Section 5.32 Calculation of $\vartheta$; Section 5.4 Elliptic and mechanical constants; Section 5.5 Discussion of elliptic arguments; Section 5.51 Discussion of functions $\Phi, \Phi_1, \Phi', \Phi_1'$;
Section 5.6 Discussion of obtained formulae; Section 5.7 The general interpretation of motion; Section 5.8 Special cases of motion.

\medskip

Chapter 6, The special cases of motion is the last one. It consists of: Section 6.1 Constants $s_0, n_0$ and  $x'$. The characteristic curve of degree 3; Section 6.2 Approximate calculation of motion; Section 6.3
Regular precession. Perturbation of motion; Section 6.4 Pseudo-regular precession; Section 6.5 Rolling of a ball over a sphere; Section 6.6 Stationary motion. Perturbed motion; Section 6.7 Remarkable trajectories.

\medskip

The first three chapters are more general and introductory. The second part, consisting of chapters 4-6 is more special, contains the original solution of the posed problem and a detailed analysis of the obtained solution. This second part occupies the major part of the text.

\subsection{Rolling without slipping of a body over a surface in the Neumann variables} \label{Nvar}

Chapter 1 introduces very convenient coordinates of C. Neumann, \cite{Neumann1899} in the study of rolling of one surface over another. Suppose that a body $\mathrm T$ bounded by its surface $\mathrm S$ is rolling over a surface $\mathrm S_1$.
Let $O_1\mathbf x_1\mathbf y_1\mathbf z_1$ be the coordinate frame fixed in the space and let $O\mathbf x\mathbf y\mathbf z$ be the frame attached to the body with the same orientation.\footnote{As positive orientation Demchenko uses nowadays negative orientation. This is the reason the signs in several equations differs from the signs we used to have.}
Let
\begin{eqnarray}
&& x_1=x_1(u_1,v_1), \qquad y_1=y_1(u_1,v_1), \qquad z_1=z_1(u_1,v_1),\label{parametrizacija1}\\
&& x=x(u,v), \qquad \qquad y=y(u,v), \qquad \qquad z=z(u,v),\label{parametrizacija}
\end{eqnarray}
be the parameterizations of $\mathrm S_1$ and $\mathrm S$ in the corresponding coordinates of the given frames, where
$u_1, v_1$ are the Gauss coordinates defined along the principle
curvature lines of $\mathrm S_1$ and similarly $u, v$ are the Gauss coordinates of $\mathrm S$.

\medskip

Let $M$ be their point of contact and $\mathbf n_1$ and $\mathbf n$ the unit vectors normal to $\mathrm S_1$ and $\mathrm S$ respectively, such that the frames $M\mathbf u_1\mathbf v_1\mathbf n_1$ and $M\mathbf u\mathbf v\mathbf n$ are positively oriented. Here $\mathbf u_1, \mathbf v_1, \mathbf u, \mathbf v$ are unit tangent vectors to the point of contact $M$ of the coordinate lines $u_1,v_1$ of $\mathrm S_1$ and the coordinate lines $u,v$ of $\mathrm S$, respectively.
Then the \emph{Neumann coordinates} of $\mathrm T$ are: the Gauss coordinates $u, v$, $u_1, v_1$ at $M$ and the angle $\vartheta$ between $\mathbf v_1$ and $\mathbf u$.

\medskip

Given moving frames $O\mathbf x\mathbf y\mathbf z$, $M\mathbf u_1\mathbf v_1\mathbf n_1$, and $M\mathbf u\mathbf v\mathbf n$, the following angular velocities are defined:
\begin{eqnarray*}
& \omega\colon  & \text{of $O\mathbf x\mathbf y\mathbf z$ with respect to $O_1\mathbf x_1\mathbf y_1\mathbf z_1$, i.e., the angular velocity of the body $\mathrm T$},\\
&  \omega_1 \colon & \text{of $O\mathbf x\mathbf y\mathbf z$ with respect to $M\mathbf u\mathbf v\mathbf n$},\\
&\omega_2  \colon &  \text{of $M\mathbf u\mathbf v\mathbf n$ with respect to $M\mathbf u_1\mathbf v_1\mathbf n_1$},\\
& \omega_3  \colon & \text{of $M\mathbf u_1\mathbf v_1\mathbf n_1$ with respect to $O_1\mathbf x_1\mathbf y_1\mathbf z_1$},
\end{eqnarray*}
related by $\omega=\omega_1+\omega_2+\omega_3$.

\medskip

The nonholonomic constraint that the body $\mathrm T$ rolls without slipping over the surface $\mathrm S_1$ is usually given by the condition that the current point of contact $M$ considered
in rest in the moving frame has also zero velocity in the space frame.\footnote{In the usual vector notation with the standard orientation,
the no-slipping condition is given in the form which is not used in the dissertation:
\begin{equation}\label{uobicajeno}
\omega \times \overrightarrow{OM}+\frac{d}{dt}\overrightarrow{O_1O}=0.
\end{equation}}
Here, the velocity of the point $M$ is different from zero. Namely, Demchenko consider the point of contact $M$ as a function of time in coordinates $u_1(t),v_1(t)$,
as a "trace" of the body $\mathrm T$ over $\mathrm S_1$. The corresponding "trace" on $\mathrm S$ is given by functions $u(t),v(t)$. Let $\mathfrak v^1$, $\mathfrak v$  be the vectors of absolute and relative velocities of $M$, i.,e,
the time derivatives of \eqref{parametrizacija1} and \eqref{parametrizacija}, respectively.
Then they are related by the expression $\mathfrak v^1=\mathfrak v+\mathfrak m$\footnote{The velocity $\mathfrak m$ is usually
expressed as the left hand side od \eqref{uobicajeno}.}.
The condition that
the body $\mathrm T$ rolls without slipping over the surface $\mathrm S_1$ is then given by
\begin{equation}\label{eq:nosliding1}
\mathfrak v^1=\mathfrak v.
\end{equation}
The condition \eqref{eq:nosliding1} is equivalent to \eqref{uobicajeno}.

\medskip

Denote by $s, \tau, n, \mathfrak m_\mathbf u,\mathfrak m_\mathbf v,\mathfrak m_\mathbf n$ the projections of the angular velocity $\omega$ and $\mathfrak m$ on the axes $\mathbf u, \mathbf v, \mathbf n$. It is clear that $\mathfrak m_\mathbf n=0$.
In Section 1.2, the formulae ((8-10)) are derived which express $s, \tau, n, \mathfrak m_\mathbf u, \mathfrak m_\mathbf v$ as homogenous linear functions of the time derivatives of the Neumann coordinates $\dot u, \dot v, \dot u_1, \dot v_1, \dot \vartheta$, and vice-versa. For further reference we will provide the equations:
\begin{eqnarray}
&& s=-\frac{D''}{\sqrt{G}}\dot v- \frac{D''_1}{\sqrt{G_1}}\dot v_1\sin\vartheta-\frac{D_1}{\sqrt{E_1}}\dot u_1\cos\vartheta, \nonumber\\
&& \tau=\frac{D}{\sqrt{E}}\dot u- \frac{D_1}{\sqrt{E_1}}\dot u_1\sin\vartheta+\frac{D''_1}{\sqrt{G_1}}\dot v_1\cos\vartheta,\nonumber\\
&& n=-\dot \vartheta + \frac{1}{2\sqrt{E_1G_1}}\big(\frac{\partial E_1}{\partial v_1}\dot u_1- \frac{\partial G_1}{\partial u_1}\dot v_1\big)+\frac{1}{2\sqrt{EG}}\big(\frac{\partial E}{\partial v}\dot u- \frac{\partial G}{\partial u}\dot v\big),\label{eq:n}
\end{eqnarray}
for which Demchenko refers to \cite{Vor1911}. Here and further $E, F, G$, $D, D', D''$ are the coefficients of the first and the second fundamental forms of $\mathrm S$. Similarly $E_1, F_1, G_1$, $D_1, D_1', D_1"$ are the coefficients of the first and the second fundamental forms of $\mathrm S_1$. The choice of the Gauss coordinates gives $F=0, D'=0$ and similarly $F_1=0, D_1'=0$.

\medskip

The condition \eqref{eq:nosliding1} gives the differential constraints $\mathfrak m_\mathbf u=0$ and $\mathfrak m_\mathbf v=0$ for rolling without sliding as expressed in (1) of Section 1.7:
\begin{equation}\label{eq:nosliding2}
\sqrt{E_1}\dot u_1= - \sqrt{E} \dot u \sin \vartheta + \sqrt{G} \dot v \cos \vartheta, \quad \sqrt{G_1}\dot v_1= \sqrt{E} \dot u \cos \vartheta + \sqrt{G} \dot v \sin \vartheta.
\end{equation}

The equations of motion of a rolling without slipping of the rigid body $\mathrm{T}$ over the surface $\mathrm S_1$ are derived in a two different equivalent ways. In Chapter 2 they are derived by using general lows of mechanics, while in Chapter 3 they are derived by using the Voronec principle.

\medskip

Let $\mathbf M$ be the mass of the body and $\mathbf w$ the velocity of the point $O$.
It is assumed that the mass center of the body $\mathrm T$ is at the point $O$ and that axes $O\mathbf x$, $O\mathbf y$, and $O\mathbf z$ are the principal axes of body. Let $p,q,r$ be the components of the angular velocity $\omega$ and $A, B, C$ be the components
of the inertia tensor, and  $\mathbf w_\mathbf x,\mathbf w_\mathbf y,\mathbf w_\mathbf z$ be the components of the velocity $\mathbf w$ in the moving frame $O\mathbf x \mathbf y\mathbf z$.
Then the kinetic energy of the body is given by
\begin{equation}\label{kineticka}
\mathbf T=\frac{\mathbf M}2\big(\mathbf w_\mathbf x^2+\mathbf w_\mathbf y^2+\mathbf w_\mathbf z^2\big) +\frac12\big(Ap^2+Bq^q+Cr^2\big).
\end{equation}
Denote the momentum of the body $\mathrm T$ as $\mathfrak M$ and the angular momentum  with respect to the point $M$ as $\mathbf G^{(M)}$.  We have general laws of mechanics
written in the fixed reference frame $O_1\mathbf x_1\mathbf y_1\mathbf z_1$:
\begin{equation}\label{eq:lawmoments}
\dot{\mathfrak M}=\mathbf F, \qquad \dot{\mathbf G}^{(M)}+ [\mathfrak v^1, \mathfrak M]=\mathbf L^{(M)},
\end{equation}
where $[\cdot,\cdot]$ is the vector product\footnote{Here, since $[\mathbf x,\mathbf y]=\mathbf z$, the sign differs from the usual one.}, $\mathbf F$ is the sum of all forces, $\mathbf L^{(M)}$ is the torque of all forces applied to the body
$\mathrm T$ with respect to $M$. Note that the forces of reactions of constraints do not induce torque with respect to $M$.

\medskip

From the equations \eqref{eq:lawmoments} written in the moving frames, \eqref{eq:nosliding2} and \eqref{eq:n},
by using quite interesting manipulations with different projections of angular velocities $\omega,\omega_1,\omega_2,\omega_3$ and derivations of the
kinetic energy $\mathbf T$ and the kinetic energy written in terms of $\mathfrak m_\mathbf u,\mathfrak w_\mathbf v,\mathfrak m_\mathbf n, s, \tau,n$,
\begin{equation}\label{kineticka*}
\bar{\mathbf T}(\mathfrak m_\mathbf u,\mathfrak w_\mathbf v,\mathfrak m_\mathbf n, s, \tau,n)=\mathbf T(\mathbf w_\mathbf x,\mathbf w_\mathbf y,\mathbf w_\mathbf z,p,q,r),
\end{equation}
the system of eight differential equations in eight unknown functions of time: $u, v, u_1, v_1, \vartheta, s, \tau, n$
(or, equivalently, $u,v,\vartheta, u_1,v_1,\dot u,\dot v,\dot\vartheta$) is derived.

\subsection{The Voronec principle}\label{VorPr}

In Chapter 3, Section 3.1, Demchenko recall on the derivation of the Voronec principle for nonholonomic systems \cite{Vor1902}.
Consider the nonholonomic system with the kinetic energy $T=T(t,q_s,\dot q_s)$ ($s=1,\dots,n+k$), the generalized forces $Q_s$ that correspond to
coordinates $q_s$, and the time-dependent nonhomogeneous nonholonomic constraints
\begin{equation}\label{v1}
\dot q_{n+\nu}=\sum_{i=1}^n a_{\nu i} (q,t) \dot q_i + a_\nu (q,t) \qquad (\nu=1,2,\dots,k).
\end{equation}

Let $\Theta$ be the kinetic energy $T$ after imposing the constraints \eqref{v1} and let $K_{\nu}$ be the partial derivative of the kinetic energy $T$
with respect to $\dot q_\nu$ restricted to the constrained subspace defined by \eqref{v1}:
\begin{eqnarray}
&& \Theta(t,q_1,\dots,q_{n+k},\dot q_1,\dots,\dot q_n)=T(t,q_1,\dots,q_{n+k},\dot q_1,\dots,\dot q_{n+k}), \label{v2} \\
&&  K_\nu(t,q_1,\dots,q_{n+k},\dot q_1,\dots,\dot q_n)=\frac{\partial T}{\partial \dot q_\nu}(t,q_1,\dots,q_{n+k},\dot q_1,\dots,\dot q_{n+k}) \quad (\nu=1,\dots,k). \label{v3}
\end{eqnarray}

Based of the Lagrange-d'Alembert principle, following Voronec \cite{Vor1902}, the equations of a motion of the given noholonomic system are derived in the form
without Lagrange multipliers:
\begin{equation}
\frac{d}{dt}\frac{\partial\Theta}{\partial \dot q_i}=\frac{\partial\Theta}{\partial q_i} + Q_i+\sum_{\nu=1}^k a_{\nu i}\big(\frac{\partial\Theta}{\partial q_{n+\nu}}+Q_{n+\nu}\big)+
\sum_{\nu=1}^k K_\nu\big(\sum_{j=1}^n A_{ij}^{(\nu)} \dot q_j+A_j^{(\nu)}\big) \quad (i=1,\dots,n), \label{v4}
\end{equation}
where the components $A_{ij}^{(\nu)}$ and $A_i^{(\nu)}$ are functions of the time $t$ and the coordinates $q_1,\dots,q_{n+k}$ given by
\begin{eqnarray*}
&& A_{ij}^{(\nu)}=  \big(\frac{\partial a_{\nu i}}{\partial q_j}+\sum_{\mu=1}^k a_{\mu j}\frac{\partial a_{\nu i}}{\partial q_{n+\mu}}\big)
                   -\big(\frac{\partial a_{\nu j}}{\partial q_i}+\sum_{\mu=1}^k a_{\mu i}\frac{\partial a_{\nu j}}{\partial q_{n+\mu}}\big), \\
&&  A_i^{(\nu)}= \big(\frac{\partial a_{\nu i}}{\partial t}+\sum_{\mu=1}^k a_{\mu }\frac{\partial a_{\nu i}}{\partial q_{n+\mu}}\big)
                   -\big(\frac{\partial a_{\nu}}{\partial q_i}+\sum_{\mu=1}^k a_{\mu i}\frac{\partial a_{\nu}}{\partial q_{n+\mu}}\big).
\end{eqnarray*}

It is interesting that the equations can be written in a compact form by using a formal expression similar to the Hamiltonian principle of least action:
\begin{equation}\label{v5}
\int_{t_1}^{t_2}\Big[\delta\Theta+\sum_{i=1}^{n+k} Q_i \delta q_i+\sum_{\nu=1}^k K_\nu\big(\frac{d}{dt}\delta q_{n+\nu}-\delta\dot q_{n+\nu}\big)  \Big]\cdot dt=0,
\end{equation}
where virtual displacements $\delta q_1,\dots,\delta q_{n}$ are arbitrary and equal to zero at the endpoints of a trajectory $q(t)$ (for $t=t_1$ and $t=t_2$), while $\delta q_{n+1},\dots,\delta q_{n+k}$ are determined from the homogeneous constraints
\begin{equation}\label{v6}
\delta q_{n+\nu}=\sum_{i=1}^n a_{\nu i} \delta q_i \qquad (\nu=1,2,\dots,k).
\end{equation}

The expression \eqref{v5} is referred as the \emph{Voronec principle}. Here $\frac{d}{dt}\delta q_{n+\nu}-\delta\dot q_{n+\nu}$ are calculated according to the expressions
\eqref{v1}, \eqref{v6} and using the rule:
\[
\frac{d}{dt}\delta q_{i}-\delta\dot q_{i}=0 \qquad (i=1,2,\dots,n).
\]

In the case when all considered objects do not depend on the variables $q_{n+1},\dots,q_{n+k}$, the system is known as the \emph{Chaplygin system} and the equations \eqref{v4}
as the \emph{Chaplygin equations}. This is the reason Bilimovi\'c used the notion Chaplygin-Voronec equations (see \cite{Andj}).

\subsection{The Voronec principle and rolling  of a body over a surface }\label{body}

In  Sections 3.2,  following  \cite{Vor1912}, Demchenko applied the Voronec principle to the above problem of rolling without slipping of a body over a surface. The nonholonomic constraints are
given by \eqref{eq:nosliding2}. One can choose $\dot u_1$ and $\dot v_1$ as a dependent velocities. The corresponding generalized impulses $K_1$ and $K_2$ are defined as
\[
K_1(u,  v, \vartheta, u_1, v_1, \dot u, \dot v, \dot \vartheta)=\frac{\partial{T}}{\partial \dot u_1},\quad K_2(u,  v, \vartheta, u_1, v_1, \dot u, \dot v, \dot \vartheta)=\frac{\partial{T}}{\partial \dot v_1}
\]
where $T={T}(u,  v, \vartheta, u_1, v_1, \dot u, \dot v, \dot \vartheta, \dot u_1, \dot v_1)$ is the kinetic energy \eqref{kineticka} in the Neumann variables and the constraints \eqref{eq:nosliding2} are imposed after the taking of partial derivatives.

\medskip

From now on, $\Theta$ denotes the kinetic energy \eqref{kineticka*} as a function of angular velocities $s, \tau, n$ taking into account the constraints $\mathfrak m_\mathbf u=0,\mathfrak w_\mathbf v=0,\mathfrak m_\mathbf n=0$, while
$\bar\Theta(u,  v, \vartheta, u_1, v_1,\dot u,\dot v,\dot \vartheta)$ denotes the kinetic energy $T(u,  v, \vartheta, u_1, v_1, \dot u, \dot v, \dot \vartheta, \dot u_1, \dot v_1)$ after imposing the constraints \eqref{eq:nosliding2}.
Then
\[
\begin{aligned}
K_1&=M\sqrt{E_1}\big[(\epsilon s-\rho_u n)\cos \vartheta+(\epsilon\tau-\rho_v n)\sin\vartheta\big]\\
&+\frac{1}{2\sqrt{E_1G_1}}\frac{\partial E_1}{\partial v_1}\frac{\partial\Theta}{\partial n}-\frac{D_1}{\sqrt{E_1}}\Big(\frac{\partial \Theta}{\partial s}\cos\vartheta+\frac{\partial \Theta}{\partial \tau}\sin\vartheta\Big),\\
K_2&=M\sqrt{G_1}\big[(\epsilon s-\rho_u n)\sin \vartheta-(\epsilon\tau-\rho_v n)\cos\vartheta\big]\\
&+\frac{1}{2\sqrt{E_1G_1}}\frac{\partial G_1}{\partial u_1}\frac{\partial\Theta}{\partial n}+\frac {D''_1}{\sqrt{G_1}}\Big(\frac{\partial \Theta}{\partial \tau}\cos\vartheta-\frac{\partial \Theta}{\partial s}\sin\vartheta\Big).
\end{aligned}
\]
 Here $\rho_u, \rho_v,\epsilon$ are the coordinates of  $\overrightarrow{OM}$ in the coordinate system $M\mathbf u\mathbf v\mathbf n$.
In these expressions, $s, \tau, n$ should be expressed as functions of $\dot u, \dot v, \dot \vartheta$ by using \eqref{eq:n} and the constraints \eqref{eq:nosliding2}.

\medskip

Having in mind that virtual displacements satisfy
\begin{equation}\label{varijacije}
\begin{aligned}
\sqrt{E_1}\delta u_1 &=-\sqrt{E}\delta u\sin\vartheta+\sqrt{G}\delta v\cos\vartheta,\\
\sqrt{G_1}\delta v_1 &=\sqrt{E}\delta u\cos\vartheta+\sqrt{G}\delta v\sin\vartheta,\\
\end{aligned}
\end{equation}
the coefficients in \eqref{v5} of terms that contain $K_1$ and $K_2$ have the form
$$
\begin{aligned}
\frac{1}{\sqrt{E_1}}\big[\sqrt{E}(n\delta u-\dot u n')\cos \vartheta +\sqrt{G}(n\delta v-\dot v n')\sin\vartheta\big],\\
\frac{1}{\sqrt{G_1}}\big[\sqrt{E}(n\delta u-\dot u n')\sin \vartheta -\sqrt{G}(n\delta v-\dot v n')\cos\vartheta\big],\\
\end{aligned}
$$
where
$$
n'=-\delta\vartheta+\frac{1}{2\sqrt{GE}}\Big(\frac{\partial E}{\partial v}\delta u-\frac{\partial G}{\partial u}\delta v\Big)+
\frac{1}{2\sqrt{G_1E_1}}\Big(\frac{\partial E_1}{\partial v_1}\delta u_1-\frac{\partial G_1}{\partial u_1}\delta v_1\Big).
$$
If one denotes
$$
K_1'=\frac{K_1}{\sqrt{E_1}}\cos\vartheta+\frac{K_2}{\sqrt{G_1}}\sin\vartheta,\quad K_2'=\frac{K_1}{\sqrt{E_1}}\sin\vartheta-\frac{K_2}{\sqrt{G_1}}\cos\vartheta,
$$
the Voronec principle \eqref{v5} can be written in the form
$$
\int\limits_{t_1}^{t}\big[\delta\bar{\Theta}+\delta U+K'_1\sqrt{E}(n\delta u-\dot u n') +K'_2\sqrt{G}(n\delta v-\dot v n')\big]dt.
$$
Here $U(u,v,\vartheta, u_1,v_1)$ is a \emph{force function} (negative potential energy).

\medskip

By using the expression for $n'$,  \eqref{varijacije}, and setting the terms that contain independent variations $\delta u, \delta v, \delta \vartheta$ to zero,
one gets the equations of the motion in the form
\begin{equation}\label{jednacine1}
\begin{aligned}
\frac{d}{dt}\frac{\bar{\Theta}}{\partial\dot u}-\frac{\partial(\bar{\Theta}+U)}{\partial u}&=\sqrt{E}\Big[-\frac{\partial(\bar{\Theta}+U)}{\partial u_1}\frac{\sin\vartheta}{\sqrt{E_1}}
+\frac{\partial(\bar{\Theta}+U)}{\partial v_1}\frac{\cos\vartheta}{\sqrt{G_1}}-K'_1\dot\vartheta\Big]\\
&-(\Delta_2K'_1+\Delta_1K'_2)\sqrt{EG}\dot v,\\
\frac{d}{dt}\frac{\bar{\Theta}}{\partial\dot v}-\frac{\partial(\bar{\Theta}+U)}{\partial v}&=\sqrt{G}\Big[\frac{\partial(\bar{\Theta}+U)}{\partial u_1}\frac{\cos\vartheta}{\sqrt{E_1}}
+\frac{\partial(\bar{\Theta}+U)}{\partial v_1}\frac{\sin\vartheta}{\sqrt{G_1}}-K'_2\dot\vartheta\Big]\\
&+(\Delta_2K'_1+\Delta_1K'_2)\sqrt{EG}\dot u,\\
\frac{d}{dt}\frac{\bar{\Theta}}{\partial\dot \vartheta}-\frac{\partial(\bar{\Theta}+U)}{\partial \vartheta}&=K'_1\sqrt{E}\dot u+K '_2\sqrt{G}\dot v,
\end{aligned}
\end{equation}
where
\begin{eqnarray*}
&& 2\Delta_1=\frac{1}{\sqrt{G}}\frac{\partial\ln E}{\partial v}-\frac{\sin\vartheta}{\sqrt{G_1}}\frac{\partial\ln E_1}{\partial v_1}-\frac{\cos\vartheta}{\sqrt{E_1}}\frac{\partial\ln G_1}{\partial u_1}, \\
&& 2\Delta_2=\frac{1}{\sqrt{E}}\frac{\partial\ln G}{\partial u}-\frac{\sin\vartheta}{\sqrt{E_1}}\frac{\partial\ln G_1}{\partial u_1}-\frac{\cos\vartheta}{\sqrt{G_1}}\frac{\partial\ln E_1}{\partial v_1}.
\end{eqnarray*}
These three differential equations of the second order, together with two constraints  \eqref{eq:nosliding2}
give a system of eight equations with eight unknown variables
$u,v,\vartheta, u_1,v_1,\dot u,\dot v,\dot\vartheta$.

\medskip

Let us mention that on the occasion of the centennial of the seminal work of Voronec \cite{Vor1912}, Russian Journal of Nonlinear Dynamics published a Russian translation of the German original, prefaced with a short, but succinct text by A. S. Sumbatov. Sumbatov indicated that Voronec went toward his principle for about 10 years. He also listed people who successfully continued the work of Voronec: Ya. Shtaerman (1915), A. Bilimovi\'c (1916), and Yu. P. Bychkov (1965-67, 2004). Let us also mention the work of Bilimovi\'c \cite{Bilim1927}, where he indicated the advantages of the Voronec equations with respect to other approaches to nonholonomic mechanics.

\subsection{Rolling of a body with a gyroscope}\label{bodyGyr}

The next step for Demchenko is to consider rolling of a body $\mathrm{T}$ with a gyroscope inside the body (Chapter 3, Section 3.3).
He assumes
that the axis of the gyroscope coincides with one of the principal axes of the body ($O\mathbf z$) and that the mass center of the body and of the gyroscope is the point $O$. It is also assumed that the forces applied to the gyroscope do not induce torque about the axis of the gyroscope.
Thus, the gyroscope will rotate with a constant angular velocity $\tilde{\omega}$ around the axes $O\mathbf z$.

\medskip

The kinetic energy of the system body + gyroscope takes the form\footnote{In the presentation of the PhD thesis we completely followed the notation of the PhD thesis \cite{Dem1924} except in denoting the total kinetic energy and the angular momentum by $\tilde{\mathbf T}$ and ${\tilde{\mathbf{G}}}^{(M)}$, respectively.}
\begin{equation}\label{kinMod}
\tilde{\mathbf T}=\mathbf T+\frac12\tilde{C}\tilde{\omega}^2=\bar{\mathbf T}+\frac12\tilde{C}\tilde{\omega}^2,
\end{equation}
where  $\mathbf T$ and $\bar{\mathbf T}$ are given by \eqref{kineticka} and \eqref{kineticka*}, such that
$p,q,r$ are the components of the angular velocity of the body $\mathrm T$, $A,B$ are $\mathbf x$ and $\mathbf y$ components of the inertia tensor of the system body + gyroscope,
$C$ and $\tilde C$ are the moments of inertia with respect to the axis $O\mathbf z$ of the body and the gyroscope
in the body frame $O\mathbf x\mathbf y\mathbf z$, and $\mathbf M$ is the mass of the system body + gyroscope (more details are given in Chapter 4, Section 4.2, see the equation \eqref{kinMod**} given below).

\medskip

The angular momentum of the system body + gyroscope with respect to the point $M$ is $\tilde{\mathbf G}^{(M)}=\mathbf G^{(M)}+\kappa$, where $\mathbf G^{(M)}$ is the angular momentum of the system with the gyroscope when the $\mathbf z$-component $\tilde\omega$ of the angular velocity of the gyroscope is set to zero, and
$\kappa =k\mathbf z$, $k=\tilde C\tilde\omega$.
Due to the presence of the gyroscope, the second equation in \eqref{eq:lawmoments}
(also written in the fixed reference frame $O_1\mathbf x_1\mathbf y_1\mathbf z_1$) takes the form
\begin{equation}\label{eq:momentgyroscope*}
\dot{\tilde{\mathbf{G}}}^{(M)}+[\mathfrak v,\mathfrak M]=\mathbf{L}^{(M)},
\end{equation}
where $\mathfrak M$ is the momentum of the system body + gyroscope and
$\mathbf L^{(M)}$ is the torque of all forces. Here, the constraint \eqref{eq:nosliding1} is imposed.

\medskip

The projections of the total angular momentum $\tilde{\mathbf G}^{(M)}$ to the axis of $M\mathbf u\mathbf v\mathbf n$ are
$$
\frac{\partial\bar{\mathbf T}}{\partial s}+k\alpha'',\quad \frac{\partial\bar{\mathbf T}}{\partial \tau}+k\beta'',\quad \frac{\partial\bar{\mathbf T}}{\partial n}+k\gamma'',
$$
where $\alpha'',\beta'',\gamma''$ are cosines of the angles between $\mathbf z$ and $\mathbf u,\mathbf v,\mathbf n$.

\medskip

It is assumed that the forces in the system are potential and given by a force function $U$.
Let $s_1,\tau_1,n_1$ and $\mathfrak v_\mathbf u,\mathfrak v_\mathbf v,\mathfrak v_\mathbf n$ be the components of $\omega_1$ and $\mathfrak v$ in the frame $M\mathbf u\mathbf v\mathbf n$.
Using the low of changing of angular momentum \eqref{eq:momentgyroscope*} and the kinematic equations $\frac{d\alpha''}{dt}=\tau_1\gamma''-n_1\beta''$ (similar for $\beta''$ and $\gamma''$), the equations are written in the form:
$$
\begin{aligned}
\frac{d}{dt}\frac{\partial\Theta}{\partial s}+(\tau-\tau_1)\frac{\partial\Theta}{\partial n}-(n-n_1)\frac{\partial\Theta}{\partial \tau}+\mathfrak v_\mathbf v\frac{\partial\bar{\mathbf T}}{\partial\mathfrak m_\mathbf n}-\mathfrak v_\mathbf n\frac{\partial\bar{\mathbf T}}{\partial\mathfrak m_\mathbf v}&=\frac{\partial\bar{\dot{U}}}{\partial s}+k(n\beta''-\tau\gamma''),\\
\frac{d}{dt}\frac{\partial\Theta}{\partial \tau}+(n-n_1)\frac{\partial\Theta}{\partial n}-(s-s_1)\frac{\partial\Theta}{\partial n}+\mathfrak v_\mathbf n\frac{\partial\bar{\mathbf T}}{\partial\mathfrak m_\mathbf u}-\mathfrak v_\mathbf u\frac{\partial\bar{\mathbf T}}{\partial\mathfrak m_\mathbf n}&=\frac{\partial\bar{\dot{U}}}{\partial \tau}+k(s\gamma''-n\alpha''),\\
\frac{d}{dt}\frac{\partial\Theta}{\partial n}+(s-s_1)\frac{\partial\Theta}{\partial \tau}-(\tau-\tau_1)\frac{\partial\Theta}{\partial s}+\mathfrak v_\mathbf u\frac{\partial\bar{\mathbf T}}{\partial\mathfrak m_\mathbf v}-\mathfrak v_\mathbf v\frac{\partial\bar{\mathbf T}}{\partial\mathfrak m_\mathbf u}&=\frac{\partial\bar{\dot{U}}}{\partial n}+k(\tau\alpha''-s\beta''),
\end{aligned}
$$
where $\frac{\partial\bar{\dot{U}}}{\partial s}$, $\frac{\partial\bar{\dot{U}}}{\partial \tau}$, $\frac{\partial\bar{\dot{U}}}{\partial n}$ are the derivatives of $U$ along the vector fields that define quasi-velocities $s,\tau,n$.

The problem reduces to the integration of eight differential equations
on eight unknown functions of time: $u, v, \vartheta, u_1, v_1, s, \tau, n$ (or, equivalently, $u,v,\vartheta, u_1,v_1,\dot u,\dot v,\dot\vartheta$).
The explicit forms of all mentioned variables and functions in term of C. Neumann variables are given in Chapters 1 and 2.

\medskip

The problem simplifies under the additional assumptions that the surface of the body is of revolution, that the axis of the gyroscope coincides with the axis of revolution and that the central ellipsoid of inertia is an ellipsoid of revolution with the axis of revolution coinciding with the axis of the gyroscope.

\subsection{The Bobilev-Zhukovsky problem and its generalization}\label{ballGyr}

If one considers rolling \emph{over the plane} and the gyroscopic body is a ball with the mass center coinciding with its geometric center, the problem can be resolved in quadratures. There are two cases
where these quadratures are elliptic. These cases were studied by Bobilev \cite{Bob1892} and Zhukovsky
\cite{Zhuk1893}. In the Bobilev case the central ellipsoid of inertia is \emph{rotationally symmetric}
and the gyroscope axis coincides the axis of symmetry, while in the Zhukovsky case the additional condition is that
\emph{the moment of the ball with  respect to the axis of the gyroscope is equal to the sum of the moments of the system ball + gyroscope with respect to the axes orthogonal to the axis of the gyroscope}.
In Chapter 4 Demchenko used the same condition as Zhukovsky and considered rolling the gyroscopic ball \emph{over a sphere}.

\medskip

Let a ball of radius $R_2$ with rotational inertia ellipsoid rolls without sliding over a fixed sphere of radius $R_1$.
The ball contains a gyroscope which axis is fixed with respect to the ball and coincides the axis of symmetry of the inertia ellipsoid of the ball.
It is assumed that the mass center of the moving ball, the mass center of the gyroscope and the geometric center of the moving ball are at the origin of the moving frame $O\mathbf x\mathbf y\mathbf z$ fixed to the ball such that $O\mathbf z$ is the axis of gyroscope.
By  $A_1, A_1, C_1$, $A_2, A_2, C_2$ are denoted the principal central moments of inertia of the ball  and  the gyroscope with respect to the frames $O\mathbf x\mathbf y\mathbf z$ attached to the ball and the frame $O\xi\eta\zeta$ rigidly connected to the gyroscope, such that $O\zeta=O\mathbf z$.
In the given notation, the \emph{Zhukovsky condition} reads
\begin{equation}\label{uslovZh}
C_1=A_1+A_2.
\end{equation}

The complements of a latitude and a longitude of the contact point $M$ are chosen for Gauss coordinates $u,v$ on the moving ball and $u_1, v_1$  on the fixed sphere. The angle $\vartheta$ is as before the angle between $u$ and $v_1$ coordinate lines.
Let $x,y,z$ and $x_1,y_1,z_1$ are coordinates of point $M$ in the moving $O\mathbf x\mathbf y\mathbf z$ and the fixed frame $O_1\mathbf x_1\mathbf y_1\mathbf z_1$, respectively. One has
$$
\begin{aligned}
& x=R_2\sin u\cos v\quad y=R_2\sin u\sin v\quad z=R_2\cos u\\
& x_1=R_1\sin u_1\cos v_1\quad y_1=R_1\sin u_1\sin v_1\quad z_1=R_1\cos u_1.
\end{aligned}
$$
Now, the nonholonomic constraints \eqref{eq:nosliding2} read
\begin{equation}\label{4.2.2}
\begin{aligned}
& \dot{u}_1=-\mu'\dot{u}\sin\vartheta+\mu'\dot{v}\cos\vartheta\sin u,\\
& \dot{v}_1\sin{u}_1=\mu'\dot{u}\cos\vartheta+\mu'\dot{v}\sin\vartheta\sin u,
\end{aligned}
\end{equation}
where $\mu'=\frac{R_2}{R_1}$.

\medskip

Let, as above, $s,\tau, n$ denote coordinates of angular velocity of the ball in the moving reference frame $M\mathbf u\mathbf v\mathbf n$.
After substitution of the constraints \eqref{4.2.2} into \eqref{eq:n}, the expressions for $s$ and $\tau$ are simplified:
\begin{equation}\label{4.2.789}
s=\mu\sin u\dot{v},\quad \tau=-\mu\dot{u},\quad n=-\dot\vartheta-\cos u\dot v-\cos u_1{\dot v}_1
\end{equation}
where $\mu=1+\frac{R_2}{R_1}=1+\mu'$.

\medskip

Let $p_1,q_1,r_1$ denote projections of angular velocity of the ball to the axes of the frame $O\mathbf x\mathbf y\mathbf z$ and
let $p_2, q_2, r_2$ denote projection of angular velocity of the gyroscope to the axes of the frame $O\xi\eta\zeta$.
It is assumed that torque of active forces for gyroscope axis are zero. Since torque of reaction of constrains for gyroscope axis is also zero, one concludes
$$
C_2r_2=k=const \quad \text{and}\quad p_1^2+q_1^2=p_2^2+q_2^2.
$$

The kinetic energy of the system ball + gyroscope
is then given by (see \eqref{kineticka})
\begin{eqnarray}
\tilde{\mathbf T}&=& \frac12\big(A_1 p_1^2+A_1 q_1^2+C_1 r_1^2\big)+\frac12\big(A_2 p_2^2+A_2 q_2^2+C_2 r_2^2\big)+\frac12\mathbf M \mathbf w^2\nonumber \\
&=& \frac12\big(A p_1^2+A q_1^2+C r_1^2\big)+\frac12 k r_2+\frac12 \mathbf M \mathbf w^2, \quad A=A_1+A_2, \quad C=C_1, \label{kinMod**}
\end{eqnarray}
where $\mathbf{M}$ is mass of the system ball + gyroscope and
$\mathbf w$ is the velocity of the point $O$.

\medskip

Since $p_1^2+q_1^2+r^2_1=s^2+\tau^2+n^2$ and for the ball we have the identity
\[
\mathbf w^2=R_2^2(s^2+\tau^2),
\]
from the Zhukovsky condition \eqref{uslovZh}, the kinetic energy of the system expressed as a function of $s,\tau,n$ takes
the form
\begin{equation}\label{kinMod*}
\tilde{\mathbf T}=\frac12(P(s^2+\tau^2)+A n^2)+\frac12 k r_2,
\end{equation}
where $P=I+A$ and $I=\mathbf{M}R_2^2$.

\medskip

The equations of a motion can be obtained from a general low of change of the angular momentum \eqref{eq:momentgyroscope*}.
Since in the considered system $[\mathfrak v, \mathfrak M]=0$ and $\mathbf L^{(M)}=0$, one concludes that the total angular momentum is constant in the fixed reference frame $O_1\mathbf x_1\mathbf y_1\mathbf z_1$:
\begin{equation}\label{GMConst}
\tilde{\mathbf{G}}^{(M)}=const.
\end{equation}
Let $\Gamma$ denotes its magnitude. One con choose the axis $O_1\mathbf z_1$, such that $\tilde{\mathbf{G}}^{(M)}=\Gamma \mathbf z_1$.
The cosines $\alpha_1'',\beta_1'',\gamma_1''$  of the angles between $\mathbf z_1$ and $\mathbf u,\mathbf v,\mathbf n$ in the Neumann variables are
\[
\alpha_1''=\sin u_1\sin\vartheta, \quad \beta_1''=-\sin u_1\cos\vartheta, \quad \gamma_1''=-\cos u_1.
\]
Thus, the projections of $\tilde{\mathbf{G}}^{(M)}$ to the axes of the moving frame $M\mathbf u\mathbf v\mathbf n$ are given by
\begin{equation}\label{integrali}
\begin{aligned}
&\Gamma\sin u_1\sin\vartheta=P s-k\sin u\\
-&\Gamma\sin u_1\cos\vartheta=P\tau\\
-&\Gamma\cos u_1=A n+k\cos u\\
\end{aligned}
\end{equation}
It is interesting to mention that the equations of a motion are obtained in the form
\begin{equation}\label{4.3.131415}
\begin{aligned}
&P\dot s-\mu'A\ n \dot{u}-P\tau(n+\cos u \dot{v})=k\mu\dot u\cos u,\\
&P\dot \tau-\mu'A\ n \sin u \dot{v}+P\ s(n+\cos u \dot{v})=k(n\sin u+\mu\dot v\sin u\cos u),\\
&A\dot n=k\mu\sin u\dot{u},
\end{aligned}
\end{equation}
by derivation of \eqref{integrali} and the
kinetic energy integral \eqref{kinMod*}.

\medskip

Finally, the problem reduces to the problem of solving the system of the eight equations \eqref{4.3.131415}, \eqref{4.2.789}, and \eqref{4.2.2}
in the variables $u, v, \vartheta, u_1, v_1, s, \tau, n$.

\subsection{Solving the system in terms of elliptic functions and elliptic integrals}

Demchenko introduces new variable $x=\cos u$ and derives an elliptic equation on $x$:
\begin{equation}\label{eq:elliptic1}
\big(\frac{dx}{dt}\big)^2=X(x),
\end{equation}
where $X(x)$ is a degree four polynomial in $x$.  Namely, integrating the last of the equations \eqref{4.3.131415}, he gets
\begin{equation}\label{4.4.4}
An=-k\mu x +C_5=-k\mu(x-x_0),
\end{equation}
where $C_5=k\mu x_0$ is a constant, as well as $x_0$. In order to get $s$, he eliminates $\tau$
from the first integrals, the area integral
\begin{equation}\label{eq:areaintegral}
P^2(s^2+\tau^2)+A^2n^2 + k^2-2k(Ps\sin u -A n\cos u)=\Gamma^2,
\end{equation}
and the kinetic energy integral
\begin{equation}\label{eq:kinetic}
P(s^2+\tau^2)+An^2=2h.
\end{equation}
He gets
\begin{equation}\label{4.4.7}
b_2s\sin u = k\mu (-b_0x^2+2b_1x_0x-\bar \Gamma),
\end{equation}
where
$$
b_0=I\mu+2A, \, b_1=I\mu +A, \, b_2=2PA,
$$
and
$$
\bar \Gamma =\frac{IC_5^2+A(\Gamma^2-k^2)-2hPA}{\mu k^2},
$$
with the inequalities
$$
b_0>b_1>P=A+1,
$$
since $\mu>1$. From the kinetic energy integral \eqref{eq:kinetic}, one gets
$$
b_2^2\tau^2=-2b_2A^2n^2+2b_2Ah-b_2^2s^2.
$$
By multiplying both sides by $\sin^2u$ and by applying the formulae \eqref{4.4.4} and \eqref{4.4.7}
one gets finally
\begin{equation}\label{4.4.12}
b_2^2\tau^2\sin^2 u=\mu^2k^2X,
\end{equation}
where
\begin{equation}\label{4.4.13}
X=2b_2(h'-x+x_0)(h'+x-x_0)(1-x^2)-(-b_0x^2+2b_1x_0x-\bar \Gamma)^2,
\end{equation}
with
$$
h'=\frac{\sqrt{2hA}}{\mu k}.
$$
By substituting $\tau$ from the second of the equations \eqref{4.2.789} into \eqref{4.4.12},
one comes to \eqref{eq:elliptic1}.

By using the first integrals of energy and area, Demchenko expresses the angular velocities $s, \tau, n$ and also $u_1$ and $\vartheta$ as functions of $x$. He needs two additional elliptic integrations
\begin{equation}\label{eq:elliptic2}
dv = \frac{\phi(x)dx}{(1-x^2)\sqrt{X}}, \quad dv_1 = \frac{F(x)dx}{\theta(x)\sqrt{X}},
\end{equation}
where $\phi, F, \theta$ are quadratic polynomials in $x$. The polynomial $X$ can be presented
in the forms:
$$
X(x)=(1-x^2)\psi(x)-\phi(x)^2=a_0(x-x^I)(x-x^{II})(x-x^{III})(x-x^{IV}),
$$
where $\psi$ is also a quadratic polynomial in $x$ and $a_0$ is a negative constant.

\subsection*{Elliptic functions and addition theorems}

Using heavily and skilfully the theory of elliptic functions as presented in \cite{Hal1888}, Demchenko
inverses the integrals \eqref{eq:elliptic1} and \eqref{eq:elliptic2}. He uses Weierstrass elliptic functions, $\wp(z)$, $\zeta(z)$, and $\sigma(z)$. The basic definitions and important identities
can be found for example in \cite{Akh4}, to list a source more modern than \cite{Hal1888}. The addition theorem for elliptic functions, in particular for the Weierstrass function played important role.

\begin{theorem}[Addition theorem]
\label{th:Weierstrass.addition} The Weierstrass
function\index{Weierstrass function} satisfies the following
addition relation:
\begin{equation}\label{eq:addition1}
\wp(u+\zeta)+\wp(u)+\wp(\zeta)
 =
\frac14
 \left(
\frac{\wp'(u)-\wp'(\zeta)}{\wp(u)-\wp(\zeta)}
 \right)^2.
\end{equation}
\end{theorem}

Some other typical identities are:
\begin{equation}\label{eq:addition2}
\wp'(\zeta)\frac{\wp'(u)-\wp'(\zeta)}{\wp(u)-\wp(\zeta)}=\wp''(\zeta)-2(\wp(u)-\wp(\zeta))(\wp(u+\zeta)-\wp(\zeta))
\end{equation}
$$
\wp(u)-\wp(v)=-\frac{\sigma(u-v)\sigma(u+v)}{\sigma(u)^2\sigma(v)^2},
$$
$$
\frac {\wp'(u)}{\wp(u)-\wp(v)}=\zeta(u-v)+\zeta(u+v)-2\zeta(u).
$$
Along with addition formulae for elliptic functions, Demchenko also used the Abel theorem
for elliptic functions, stating that the sum of zeros of an elliptic functions equals the sum
of poles (modulo the lattice which defines the underlying elliptic curve).

\subsection*{Inversion of elliptic integrals}

In order to integrate the equation \eqref{eq:elliptic1}, Demchenko used an approach explained in \cite{Hal1888}, which is based on simultaneous parameterizations of the square of the polynomial $X$ of degree four in $x$ and the variable $x$ in terms of elliptic functions of the same argument $u$.
To that end, let us denote
$$
2y = \frac{\wp'(u)-\wp'(\zeta)}{\wp(u)-\wp(\zeta)}.
$$
Using the Addition Theorem \ref{th:Weierstrass.addition} and formulae \eqref{eq:elliptic1} and \eqref{eq:elliptic2}, one gets
$$
y^2-3\wp(\zeta)=(\wp(u)-\wp(\zeta))+(\wp(u+\zeta)-\wp(\zeta)),
$$
and
$$
\wp''(\zeta)-2y\wp'(\zeta)=2(\wp(u)-\wp(\zeta))(\wp(u+\zeta)-\wp(\zeta)).
$$
Let us introduce the polynomial $Y$ of degree four in $y$ as:
$$
Y=(y-3\wp(\zeta))^2+2(2y\wp'(\zeta)-\wp''(\zeta).
$$
From the above formula it follows that
$$
Y=(\wp(u+\zeta)-\wp(u))^2.
$$
Thus:
$$
\sqrt{Y}=\wp(u+\zeta)-\wp(u),
$$
and
$$
Y= y^4-6y^2\wp(\zeta)+4y\wp'(\zeta)+9\wp(\zeta)-2\wp''(\zeta).
$$
One can apply the above considerations to an arbitrary polynomial $X$ of degree four in $x$:
$$
X(x)=a_0x^4+4a_1x^3+6a_2x^2+4a_3x+a_4.
$$
To eliminate the second term with $x^3$ one substitutes the variable $x=y+h$, i.e. $x=y-a_1/a_0$.
One gets
$$
\wp(\zeta)=\frac{{a_1}^2-a_0a_2}{{a_0}^2}, \quad \wp'(\zeta)=\frac{{a_0}^2a_3-3a_0a_1a_2+2a_1^3}{a_0^3}
$$
and
$$
x=-\frac{a_1}{a_0}+\frac{1}{2}\frac{\wp'(u)-\wp'(\zeta)}{\wp(u)-\wp(\zeta)}, \quad \sqrt{X}=\sqrt{a_0}(-\wp(u+\zeta)+\wp(u)).
$$
From the Addition Theorem \ref{th:Weierstrass.addition}, it also follows that
$$
-\wp(u+\zeta)+\wp(u)=\frac{1}{2}\frac{d}{du}
\frac{\wp'(u)-\wp'(\zeta)}{\wp(u)-\wp(\zeta)},
$$
and
$$
\sqrt{X}=\sqrt{a_0}\frac{dx}{du}.
$$
Finally.
$$
\frac{u}{\sqrt{a_0}}=\int\frac {dx}{\sqrt X}.
$$

\subsection{Back to the Demchenko case}

In  general, the polynomial $X$ can have zero, two, or four real roots. The first case would not produce any real motion and Demchenko did not consider it.

\begin{figure}[h]
\centering
{\includegraphics[width=15cm]{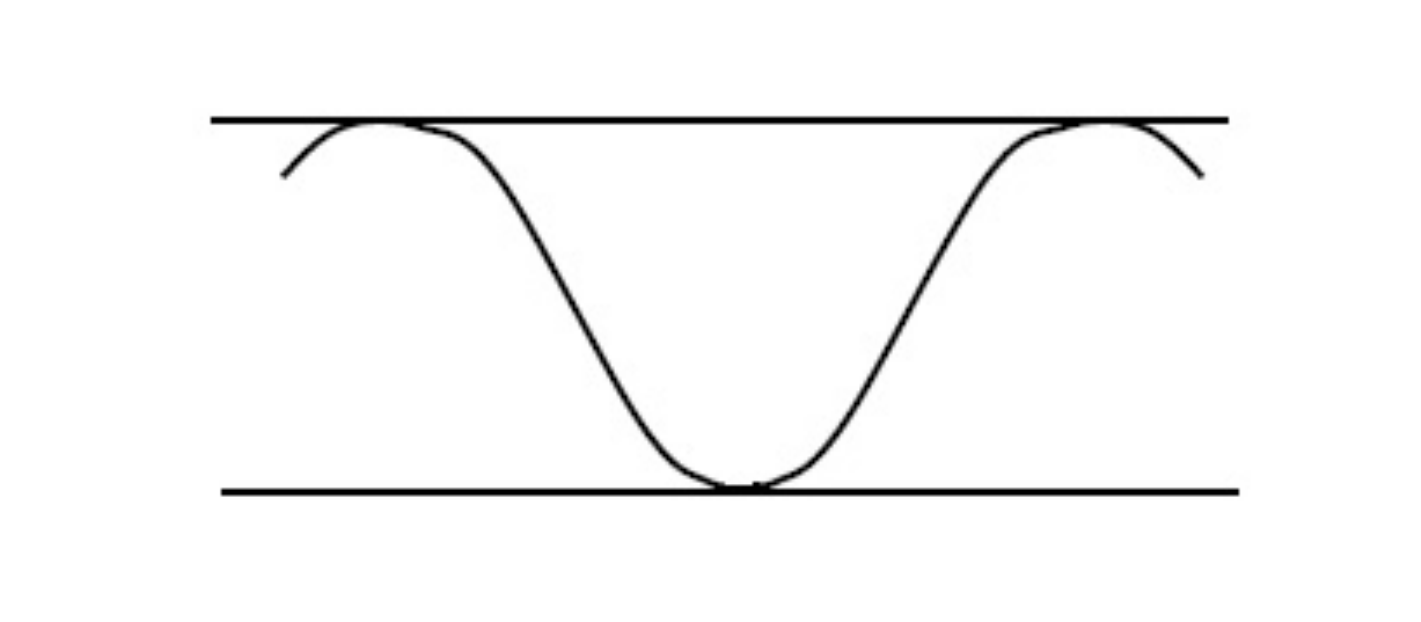}}
\caption{Demchenko: Figure 4 p. 53: the curve $A$.}\label{fig:slika4}
\end{figure}

\begin{figure}[h]
\centering
{\includegraphics[width=15cm]{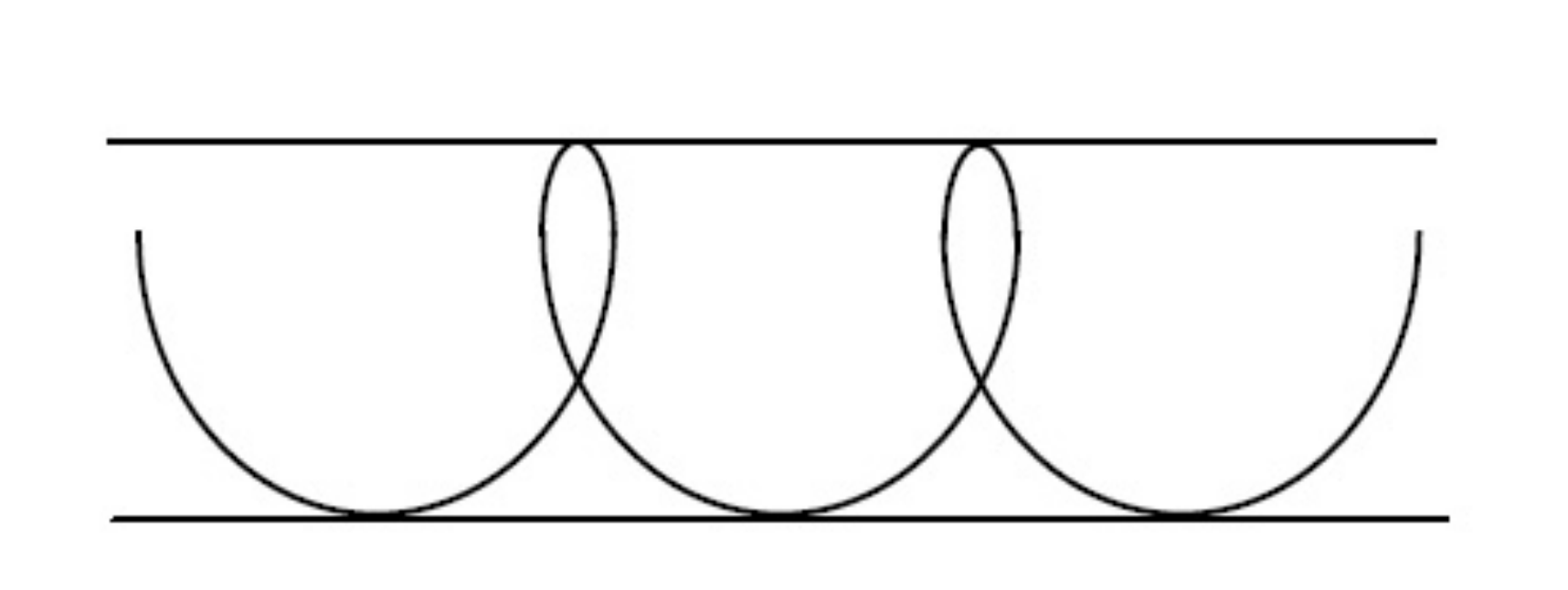}}
\caption{Figure 5, p. 53: the curve $B$.}\label{fig:slika5}
\end{figure}

\begin{figure}[h]
\centering
{\includegraphics[width=15cm]{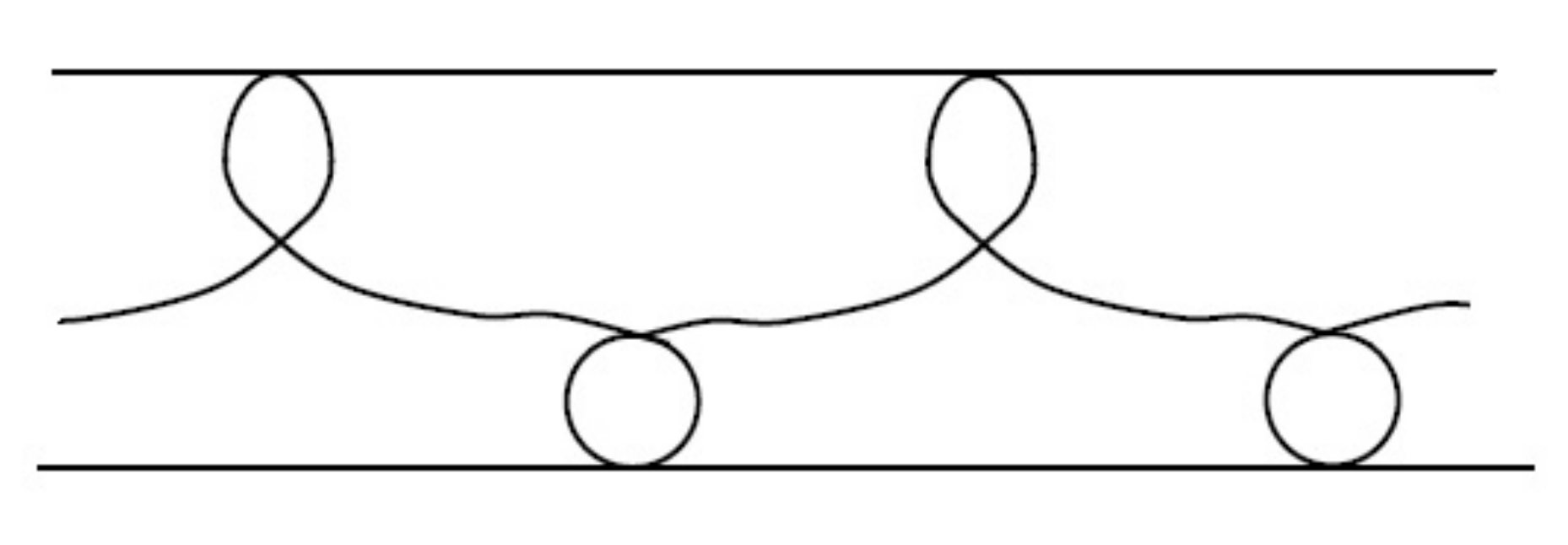}}
\caption{Demchenko: Figure 6 p. 53: the curve $C_1$.}\label{fig:slika6}
\end{figure}

\begin{figure}[h]
\centering
{\includegraphics[width=15cm]{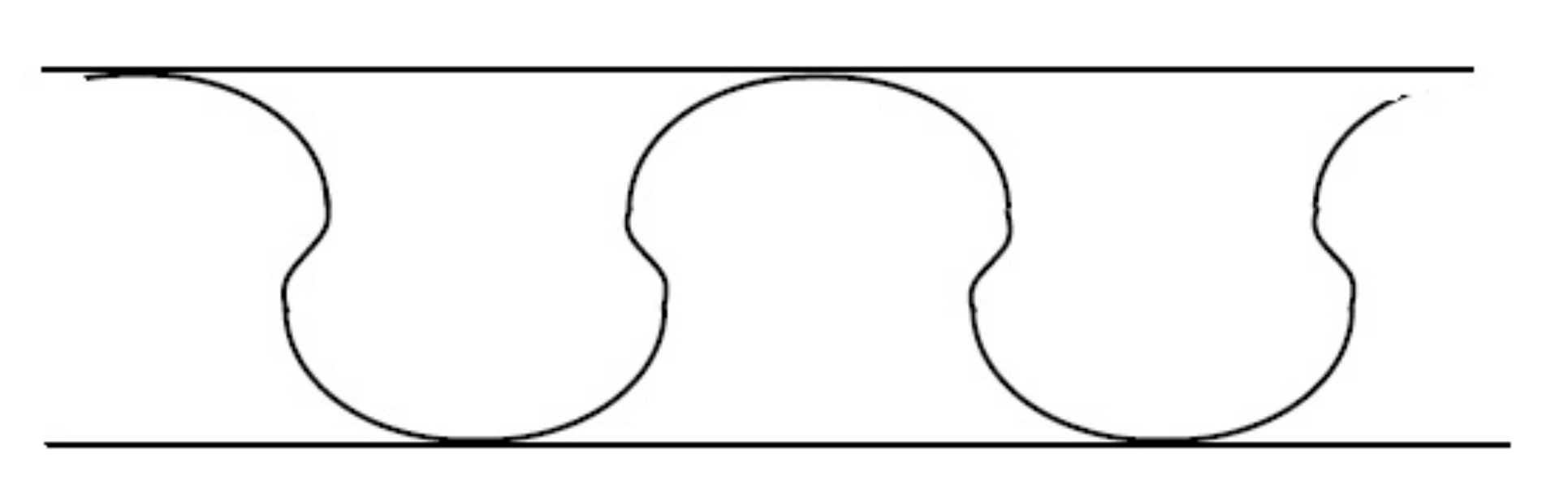}}
\caption{Demchenko: Figure 7, p. 53: the curve $C_2$.}\label{fig:slika7}
\end{figure}

\medskip
In the case of four real roots, ordered $x^{I}>x^{IV}>x^{III}>x^{II}$ the motion is possible for
$x\in (x^{I}, x^{IV})$ or $x\in (x^{III}, x^{II})$. Without loosing the generality, Demchenko works
with the first case:  $x\in (x^{I}, x^{IV})$. In the case of two real roots, he again denotes them as $x^{I}> x^{IV}$. The trajectory of the point $M$ on the mobile sphere goes between to parallels $u^{I}$ and $u^{IV}$ which it touches alternatively. The distance between two consecutive points of contact is constant. Demchenko distinguishes three cases:

\begin{itemize}
\item [1)] The polynomial $\phi(x)$ has no roots in the interval $(x^{I}, x^{IV})$.
The situation in this case is presented as curve $A$, see Figure \ref{fig:slika4}.
\item [2)] The polynomial $\phi(x)$ has one root in the interval $(x^{I}, x^{IV})$.
The situation in this case is presented as curve $B$, see Figure \ref{fig:slika5}.
\item [3)] The polynomial $\phi(x)$ has two roots in the interval $(x^{I}, x^{IV})$.
The situation in this case is presented as curves $C$, $C_1$ and $C_2$,  see  Figure \ref{fig:slika6} and Figure \ref{fig:slika7}.
\end{itemize}

The trajectories of the point $M$ on the fixed sphere are similar, where the number of roots of the polynomial $F(x)$ in the interval $(x^{I}, x^{IV})$ now discriminates cases $A$, $B$, and $C$.

\begin{figure}[h]
\centering
{\includegraphics[width=15cm]{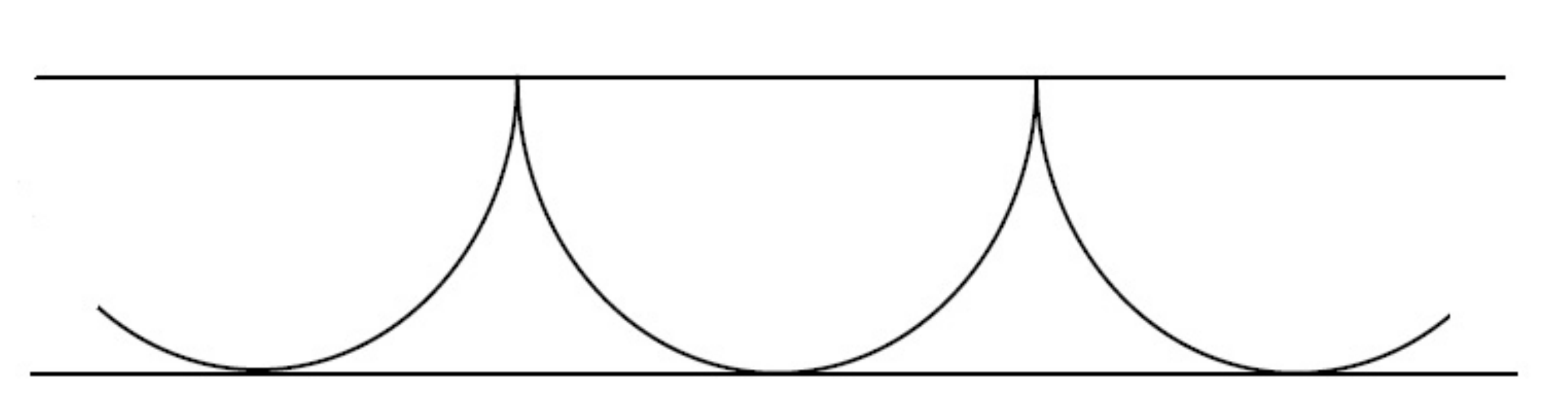}}
\caption{Demchenko: Figure 8 p. 56: the curve $D_1$.}\label{fig:slika8}
\end{figure}

\begin{figure}[h]
\centering
{\includegraphics[width=15cm]{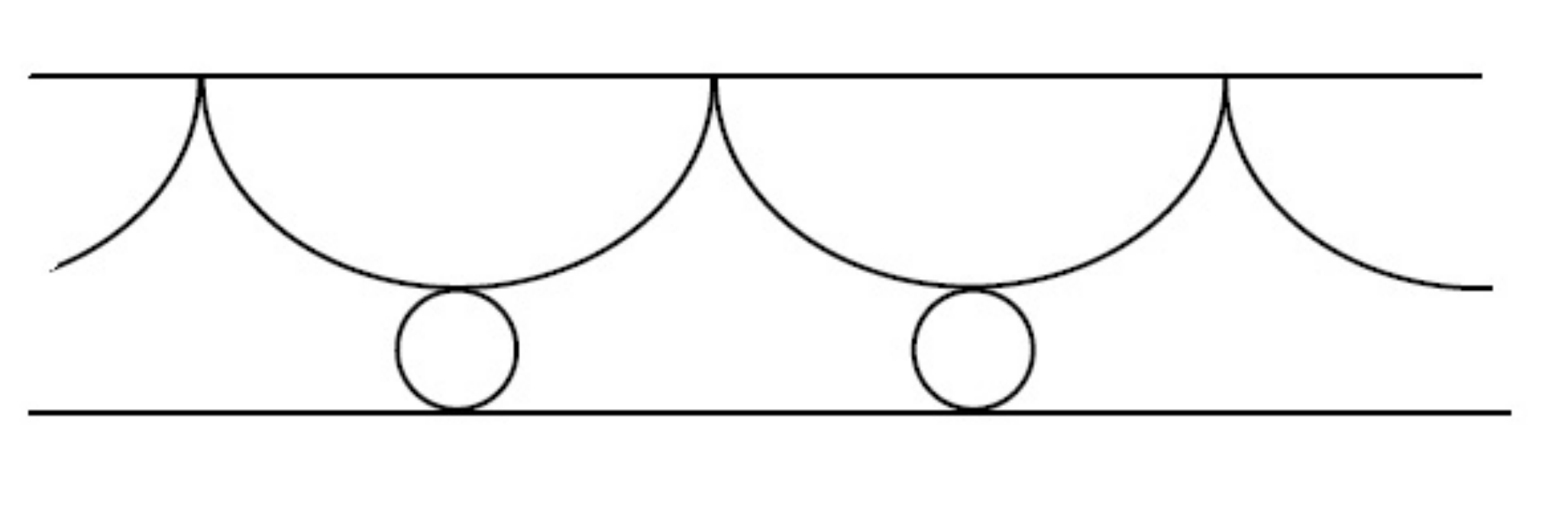}}
\caption{Demchenko: Figure 9, p. 56: the curve $D_2$.}\label{fig:slika9}
\end{figure}

\begin{figure}[h]
\centering
{\includegraphics[width=15cm]{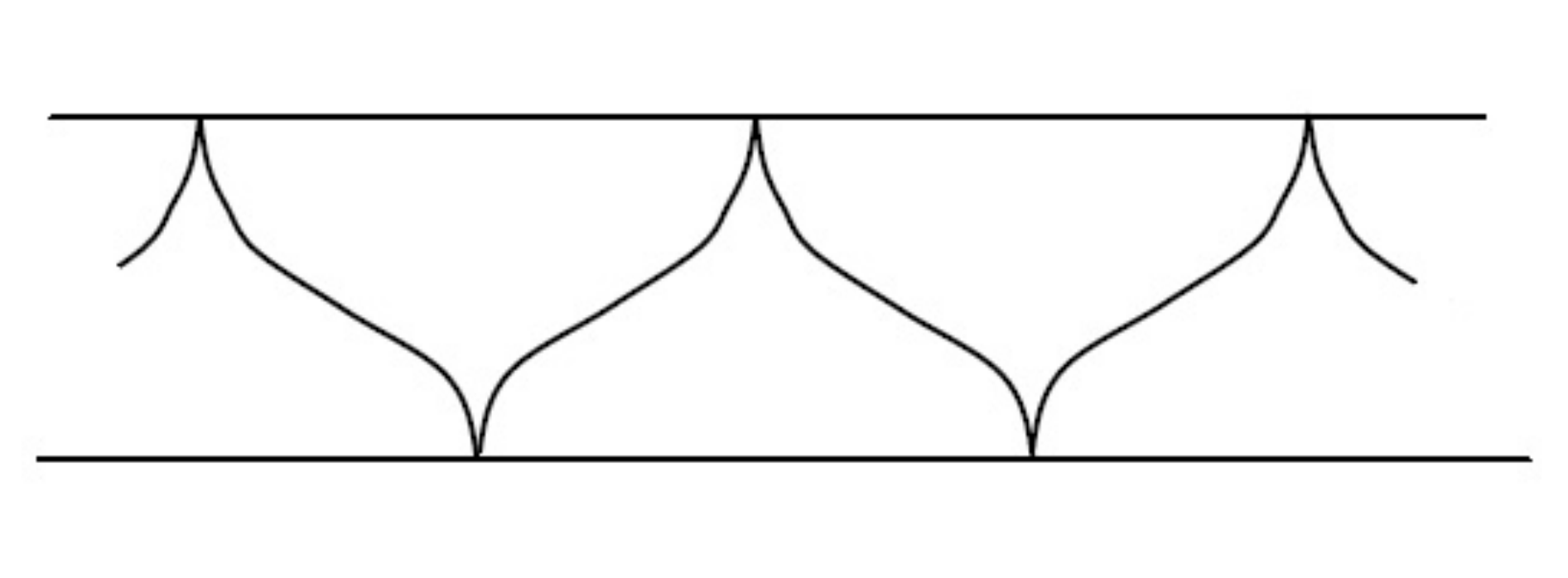}}
\caption{Demchenko: Figure 10 p. 57: the curve $D_3$.}\label{fig:slika10}
\end{figure}

\begin{figure}[h]
\centering
{\includegraphics[width=10cm]{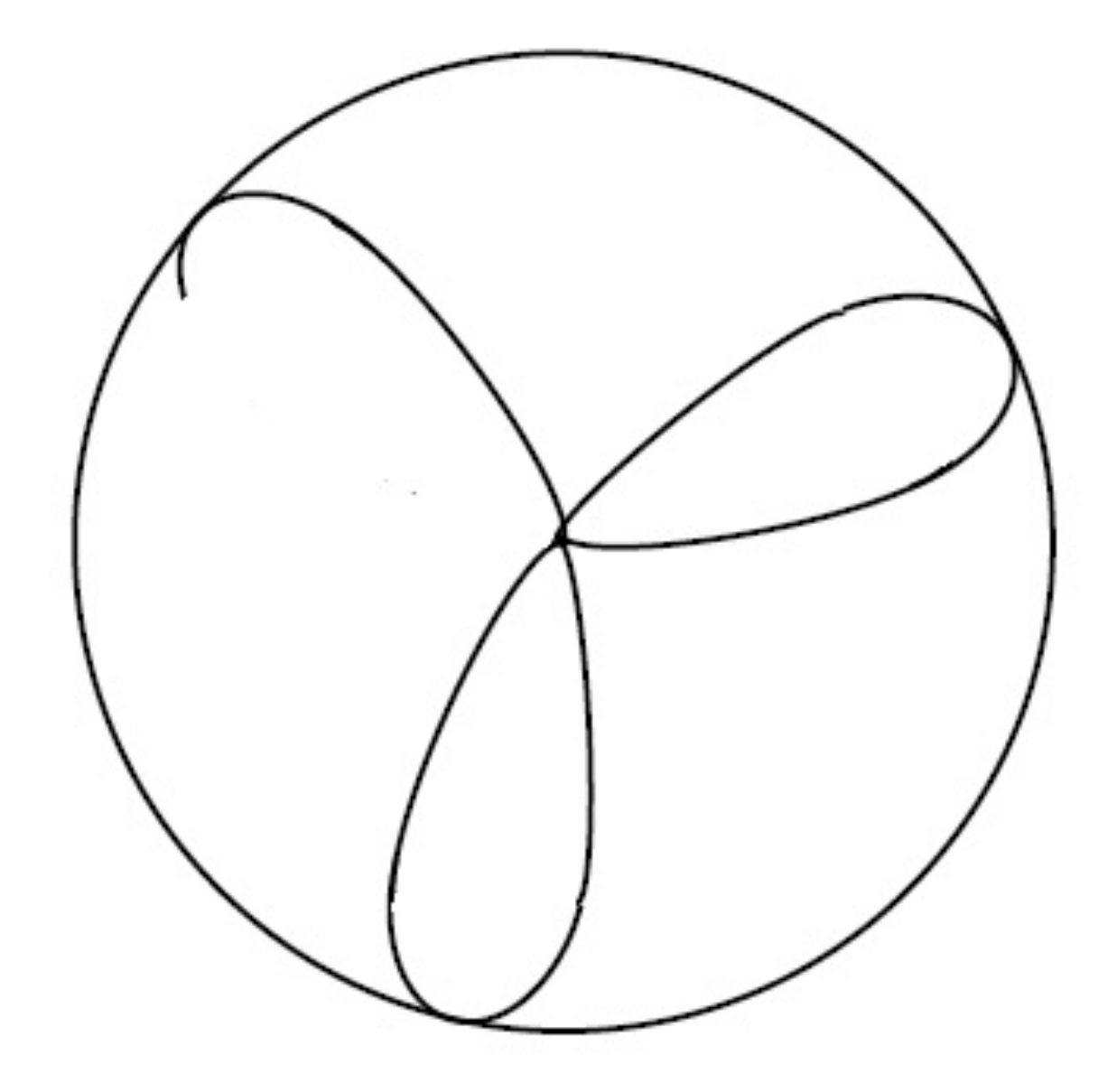}}
\caption{Demchenko: Figure 11, p. 57: the curve $E_1$.}\label{fig:slika11}
\end{figure}

\begin{figure}[h]
\centering
{\includegraphics[width=10cm]{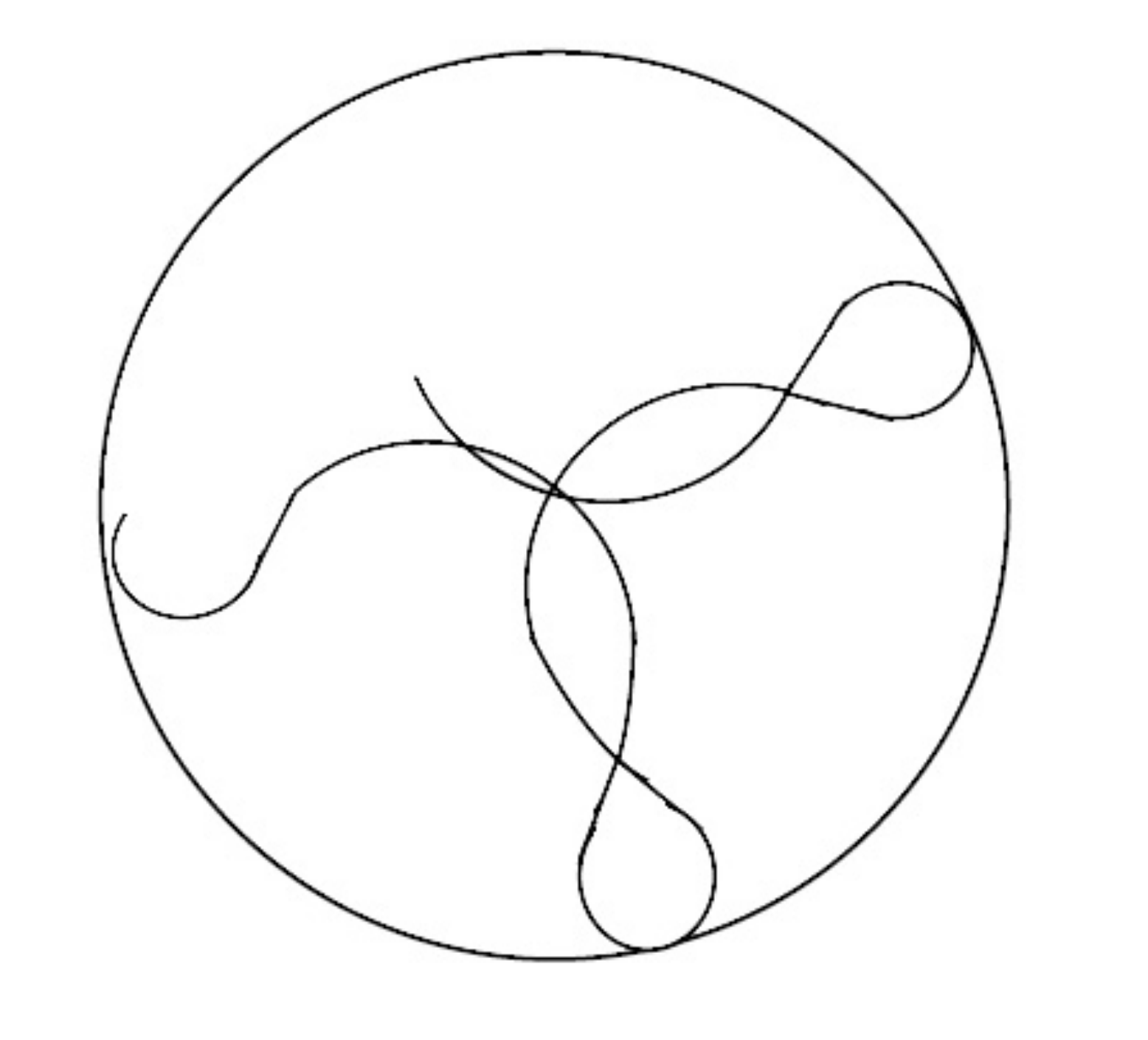}}
\caption{Demchenko: Figure 12, p. 58: the curve $E_2$.}\label{fig:slika12}
\end{figure}

There are special cases of curves if $x^{I}$ or $x^{IV}$ coincides with one of the roots of the polynomial
$\psi(x)$ or are equal to $\pm 1$. If  $x^{I}$ or $x^{IV}$ coincides with one of the roots of the polynomial
$\psi(x)$, then the curves on the movable and fixed sphere have the form $D$: $D_1$ see Figure \ref{fig:slika8}, $D_2$ see Figure \ref{fig:slika9}
and $D_3$ see Figure \ref{fig:slika10}. If, however $x^{I}=1$ or $x^{IV}=-1$ the curves representing the motion of the point $M$ on the movable sphere are presented as $E_1$ see Figure \ref{fig:slika11} and $E_2$ see  Figure \ref{fig:slika12}. In these cases the curves on the fixed sphere do not posses singularities.




\subsection{Special Solutions}

The last chapter is devoted to particular cases and particular solutions. These considerations reduce to more elementary situations than the general ones or to approximate formulae. Demchenko heavily used the capital four volume treatise of dynamics of top by Felix Klein and Arnold Sommerfeld \cite{KS}. Demchenko establishes four classes of particular motion. In each class, he also resolves the issue of stability.

The classes are:

\begin{itemize}
\item[1)] Regular precessions: they are possible. The curves which $M$ describes on both fixed and movable spheres coincide with parallels. The precession is stable if
    $$
    \frac{\partial^2 X}{\partial x^{I2}}<0,
    $$
    and instable if
    $$
    \frac{\partial^2 X}{\partial x^{I2}}>0.
    $$
In the cases
$$
    \frac{\partial^2 X}{\partial x^{I2}}=0,
    $$
    the stable situation corresponds to
    $$
    \frac{\partial^3 X}{\partial x^{I3}}=0,
    $$
    and the unstable
    the stable situation corresponds to
    $$
    \frac{\partial^3 X}{\partial x^{I3}}\ne0.
    $$
    \item[2)] Pseudo-regular precessions are possible if the rotations of the gyroscope are much bigger than the initial rotations $s_0, \tau_0$ of the gyroscopic ball. These precessions are always stable.
        \item [3)] Stationary motions are possible. The trajectories of the point $M$ consists of a single point both on movable and fixed spheres. The angular velocity of the ball is constant and the orientation of the axis of the gyroscope is constant. Such motion is always stable.
            \item [4)] Rolling of an ordinary ball, when the gyroscope stays at rest.
\end{itemize}

The last Subchapter 6.7 is devoted to the remarkable trajectories, the notion introduced by Painlev\'e \cite{Pain} and Bilimovi\'c \cite{Bilim}. These are trajectories independent on the initial energy. A detailed analysis shows that in the dynamics of gyroscopic ball rolling without sliding over a sphere such remarkable trajectories exist. These are regular precessions when the axis of the gyroscope rests parallel to the fixed vector of the moment of the gyroscopic ball with respect to the point $M$.

\section{V. S. Zhardecki, K. P. Voronec, the Paris period and fluid mechanics}

A few months before Demchenko, another immigrant from Russia, Viachislav Sigmundovich Zhardecki (1896-1962) defended his thesis
\cite{Zhar1923} also on rigid body dynamics, and having Anton Bilimovi\' c as the advisor and Milutin Milankovi\'c as a co-signatory of the report. Zhardecki family  belongs to the Polish nobility. Viachislav knew Bilimov\'c from their Odessa time and the Novorossisk University. In Belgrade he was also influenced by Milankovi\'c, a notable geoscientist and mathematician. Thus, later on, Zhardecki shifted his interests more toward geoscience and obtained remarkable results, see for example \cite{Jar1935}, \cite{Jar1948} and \cite{Jar1958}. In 1943, during the German occupation of Belgrade,  Zhardecki refused to serve at the reformed University. As a consequence, he got retired at the age of 48.
He managed to move to Austria in 1944 and in 1946-47 he served as Acting Director of the Institute of Physics and Astronomy in Graz. During that period he got experimental confirmation of his theory of formaton of continents and moved further to the US, see \cite{BlRi}. His son, Oleg, who became a scientist himself left
interesting notes about his family and the dramatic time of their emigration from Russia, see \cite{ZharO2008}.

\medskip

After defending his thesis, V. Demchenko taught mathematics in Subotica, a city 200 km north of Belgrade. His father, Grigorij Vasilievich Demchenko was a Professor of Law School there, and served as Dean of the School in 1929-30.  Both the father and the son were delegates of the Congress of Russians from Abroad in Paris in 1926, as representatives of the Yugoslav Committee, \cite{EmRus}. The same year, Bilimovi\'c and Demchenko reacted together on the paper \cite{Tz} of the notable Bulgarian scientist Ivan Cenov  and indicated three papers of P. Voronec \cite{Vor1902}, \cite{Vor1903}, \cite{Vor1911}, three papers of Bilimovi\'c \cite{Bilimovic1914}, \cite{Bilimovic1913a} \cite{Bilim} and the doctoral thesis of Demchenko \cite{Dem1924} as relevant and being source of the results close to those presented in \cite{Tz}. Their remark was published in the Liouvile's journal editorial comment \cite{Red}.

\medskip

Around that time, Vasilie
moves to Paris.
Vasilie Demchenko, now as Basile Demtchenko, defended his second doctoral dissertation in mathematical sciences in Paris in 1928. He switched his field from nonholonimic to fluid mechanics.  The thesis was entitled "I. Sur les cavitations solitaires dans un liquid infini. II. Sur l'influence des bords sur mouvement d'un corps solide dans une liquide." It was defended on June 2nd with the committee consisting of three major French mathematicians, Paul Painlev\'e (1863-1933), Henri Villat (1879-1972), and Paul Montel (1876-1975). The thesis was completed under the direction of Painlev\'e and was dedicated to Peter Voronec, the teacher: "A mon cher et regrett\'e Maitre, Pierre Voronetz."

\medskip

Demchenko's work in Paris was associated with the group of the renowned expert in hydro and aerodynamics, D. P. Ryabushinsky (1882--1962), see \cite{Andj}. Some of notable works of Demchenko include \cite{Dem1933}, \cite{Dem1937}, and \cite{Dem11938}.

\medskip

Vasilije Demchenko was an invited speaker of the International Congresses of Mathematicians \cite{ICM} in Bologna 1928 and Zurich 1932.
His two members of the Belgrade thesis committee, Bilimovi\'c and Petrovi\'c were also invited speakers at the same Congresses. In addition, Petrovi\'c was also invited speaker in Rome 1908, Cambridge 1912, and  Toronto 1924. (Demchenko's advisor from Paris, P. Painlev\'e was plenary speaker in Heidelberg 1904.)

\medskip

There is an interesting parallel between Demchenko and Konstantin Voronec, the above mentioned son of P. Voronec.
Konstantin defended his doctoral dissertation in Belgrade in 1930 \cite{Vor1930}, having the same committee as Demchenko's thesis, Bilimovi\' c, Milankovi\'c, Petrovi\'c. The Voronec thesis was very much influenced by Demchenko's thesis. After the defence, Konstantin also moved to Paris and also switched to fluid mechanics. He also defended his second doctoral thesis in Paris \cite{Vor1934a}, \cite{Vor1934b} see \cite{EmRus}, \cite{SDj}. In his second thesis, Voronec was again influenced by Demchenko, this time by \cite{Dem1933}.


\section{Demchenko's PhD thesis and contemporary nonholonomic mechanics}

The doctoral dissertation of Tatomir Andjeli\'c  can be seen as one of the important links
between the works of Voronec, Bilimovi\'c, and Demchenko and the contemporary science, \cite{Andj}. The thesis was completed just before the second word war but was defended after the war, in 1946. That was one more example of the principle adopted by many notable Serbian scientists not to participate in the university matters during the German occupation.  Although formally Bilimovi\'c didn't serve as a committee member, Andjeli\'c made it clear that the problem was posed by Bilimovi\'c and was written under his guidance. He studied application of the Voronec principle to the problem of motion  of a nonholonomic system placed in an incompressible fluid.

\medskip

An important reference in nonholonomic mechanics after the Second World War is the monograph by Neimark and Fufaev \cite{NeFu}.
There is a whole chapter devoted to the Voronec and Chaplygin equations. Among others,  the monograph referred  to several contributions of members of the Bilimovi\'c school \cite{Andjelic1954,  Bilimovic1910,Bilimovic1913a, Bilimovic1913b, Bilimovic1914, Bilimovic1915, Veljko1964a, Veljko1964b}.

\medskip

Let us note that the Chaplygin systems have a natural geometrical framework -- the nonholonomic constraints define connections on principal bundles (see  Koiler \cite{Koi1992}). On the other hand,  Bloch, Krishnaprasad, Marsden, and  Murray \cite{BKMM}  incorporated nonholonomic systems into the geometrical framework
of the Ehresmann connections. 
It was pointed out in Bak\v sa \cite{Baksa2012} that the equations used in \cite{BKMM} are literally the same as the original Voronec equations \cite{Vor1902}.
The same year de Leon also referred to the Voronets equations in \cite{deLeon}.
Now we can say that the Voronec equations,  together with the Chaplygin equations and the
 the equations of the nonholonomic systems written in terms
of quasi-velocities, known as \emph{the Euler-Poincar\'e-Chetayev-Hamel} equations,
form the central tools in the study of nonholonomic systems (e.g., see \cite{NeFu, BKMM, EKMR2005, EK2019, Zenkov2016}).

\medskip

Consider a Lagrangian nonholonomic system $(Q,L,\mathcal D)$ where the constraints define a nonintegrable
distribution $\mathcal D$ of the tangent bundle $TQ$, i.e, the constraints are homogeneous and do not depend of time.
Further we assume that $Q$ has a structure of the fiber bundle $\pi\colon Q\to S$ over the base space $S$ and that $\mathcal D$ is transverse to the fibers
of $\pi$:
\[
T_q Q=\mathcal D_q \oplus \mathcal V_q, \qquad \mathcal V_q=\ker d\pi(q).
\]
The space $\mathcal V_q$ is called the \emph{vertical space} at $q$. The distribution $\mathcal D$ can be seen as the kernel of
a vector-valued one form $A$ on $Q$, which defines the \emph{Ehresmann connection}, that satisfies

\medskip

(i) $A_q\colon T_q Q \to\mathcal V_q$ is a linear mapping, $q\in Q$;

(ii) $A$ is a projection: $A(X_q)=X_q$, for all $X_q\in\mathcal V_q$.

\medskip

The distribution $\mathcal D$ is called \emph{the horizontal space} of the Ehresmann connection $A$.
By $X^h$ and $X^v$ we denote the horizontal and the vertical component of the vector field $X\in\mathfrak{X}(Q)$.
The \emph{curvature} $B$ of the connection $A$ is a vertical vector-valued two-form defined by
\[
B(X,Y)=-A([X^h,Y^h])
\]
In the local coordinates, we have
\begin{eqnarray*}
&& \pi\colon (q_1,\dots,q_n,q_{n+1},\dots,q_{n+k}) \longmapsto (q_1,\dots,q_n),\\
&& A=\sum_{\nu=1}^k \omega^\nu \frac{\partial}{\partial q_{n+\nu}}, \quad \omega^\nu=dq_{n+\nu}-\sum_{i=1}^n a_{\nu i} dq_i, \\
&& B=\sum_{\nu=1}^k B^\nu\frac{\partial}{\partial q_{n+\nu}}, \quad B^\nu=\sum_{1\le i<j\le n} B^\nu_{ij} dq_i\wedge dq_j,\\
&& B_{ij}^{\nu}=  - \big(\frac{\partial a_{\nu i}}{\partial q_j}+\sum_{\mu=1}^k a_{\mu j}\frac{\partial a_{\nu i}}{\partial q_{n+\mu}}\big)
                   +\big(\frac{\partial a_{\nu j}}{\partial q_i}+\sum_{\mu=1}^k a_{\mu i}\frac{\partial a_{\nu j}}{\partial q_{n+\mu}}\big),
\end{eqnarray*}
i.e., $B_{ij}^\nu=-A^{(\nu)}_{ij}$ in the Voronec equations \eqref{v4}. Also,
in the case when the generalized forces $Q_s$, $s=1,\dots,n+k$ are potential: $Q_s=-{\partial V}/{\partial q_s}$, the Voronec equations \eqref{v4}
take the form:
\begin{equation}
\frac{d}{dt}\frac{\partial L_c}{\partial \dot q_i}=\frac{\partial L_c}{\partial q_i} +\sum_{\nu=1}^k a_{\nu i}\frac{\partial L_c}{\partial q_{n+\nu}}
-\sum_{\nu=1}^k\sum_{j=1}^n\frac{\partial L}{\partial q_{n+\nu}} B_{ij}^{\nu} \dot q_j \quad (i=1,\dots,n), \label{v4*}
\end{equation}
where the Lagrangian is the difference of the kinetic and the potential energy $L(q,\dot q)=T(q,\dot q)-V(q)$, and $L_c$ is the constrained Lagrangian
$L_c=L(q,\dot q^h)=\Theta-V$.
The Voronec principle \eqref{v5} for the equations \eqref{v4*} in an invariant form can be expressed as (see \cite{BKMM}):
\begin{eqnarray}\label{v5*}
\delta L_c=\mathbb FL( B(\dot q,\delta q))
\end{eqnarray}
for all virtual displacements
\[
\delta q=\sum_{s=1}^{n+k} \delta q_s\frac{\partial}{\partial q_s}\in \mathcal D_q.
\]
Here $\delta L_c$ is the variational derivative of the constrained Lagrangian along the variation $\delta q$ and $\mathbb FL$ is the fiber derivative of $L$:
\begin{equation*}
\delta L_c=\sum_{s=1}^{n+k}\big(\frac{\partial L_c}{\partial q_s}-\frac{d}{dt}\frac{\partial L_c}{\partial \dot q_s}\big)\delta q_s,
\quad \mathbb FL( B(\dot q,\delta q))=\sum_{\nu=1}^k \frac{\partial L}{\partial \dot q_{n+\nu}}B^\nu(\dot q,\delta q),
\end{equation*}
In the case when the constraints are nonhomogeneous and time dependent \eqref{v1}, the coefficients $A^{(\nu)}_{ij}$, $A^{(\nu)}_i$ can be also
interpreted as the component of the curvature of the Ehresmann connection of the fiber bundle $\pi: Q\times \mathbb R\to S\times \mathbb R$ (see Bak\v sa \cite{Baksa2012}).

\medskip

Assume that the fibration $\pi: Q\to S$ is determined by a free action of a Lie group $G$ on $Q$ ($S=Q/G$) and that the constraint distribution $\mathcal D$ and
the Lagrangian $L=T-V$ are $G$--invariant. Then $A$ is a principal connection and
the nonholonomic system \eqref{v5*} is $G$--invariant and reduces to the tangent bundle of the base manifold $S$.
The equations take the form
\begin{equation}
\delta L_{red}=\sum_{i=1}^n\big(\frac{\partial L_{red}}{\partial x_i}-\frac{d}{dt}\frac{\partial L_{red}}{\partial
\dot x_i}\big)\delta x_i= JK(\dot x,\delta x)\quad
\mathrm{for\;all} \quad \delta x\in T_x S,
\label{ChaplyginRed}
\end{equation}
where the reduced Lagrangian
$L_{red}$ is obtained from the constrained Lagrangian $L_c$ by the identification $TS=\mathcal D/G$,
and $JK(X,Y)$ is a $(0,2)$--tensor field on the base manifold $S$,
which depends on the metric and the curvature of the connection, induced by the right hand side of  \eqref{v5*}.

\medskip

The system $(Q,L,\mathcal D,G)$ is referred to as \emph{a $G$--Chaplygin system}, as a generalization of the classical Chaplygin systems
with Abelian symmetries \cite{Chaplygin1911, NeFu, Baksa, Koi1992, CCL2002, Naranjo2020}.

\medskip

Demchenko noticed that Voronec derived his principle in order to relate the nonholonomic systems to the Hamiltonian variational principle
of least action. Obviously, Voronec and his followers were aware of the fact that the equations were not variational, or, in a modern terminology,
that they were not Hamiltonian.
However, as it was pointed out by Chaplygin \cite{Chaplygin1911} that
some systems have an invariant measure, which puts them rather close to Hamiltonian
systems.  The existence of an invariant measure for various nonholomic problems is well studied  (e.g., see
\cite{Fe1988, VeVe2, Kozlov,  FK1995, ZenkovBloch, FNM, Jov2015, FGM2015}).
A closely related problem is
the Hamiltonization of nonholonomic systems
(e.g., see \cite{Chaplygin1911, BM, BM2008, BN, BBM, BMT, BBM2013, EKMR2005, CCL2002, EKMR2005, FJ2004})
Chaplygin was also
one of the first who considered a time reparametrization in order
to transform nonholonomic systems to the Hamiltonian form
\cite{Chaplygin1911}.
In the case of integrability, the dynamics over
regular invariant $m$--dimensional tori, in the original time,
has the form
\begin{equation}
\dot\varphi_1=\omega_1/\Phi(\varphi_1,\dots,\varphi_m),\dots,
\dot\varphi_m=\omega_m/\Phi(\varphi_1,\dots,\varphi_m), \qquad \Phi>0.
\label{Jacobi2}
\end{equation}

Also, after \cite{Ch1}, one of the most famous
solvable problems in nonholonomic mechanics, describing the
rolling without slipping of a balanced ball over a horizontal
surface, is referred as the \emph{Chaplygin ball}, see \cite{Kozlov2, AKN, BoMa, BM2008}.
On the other hand, the rolling without slipping of the Chaplygin ball over a sphere generically is not integrable. The only known integrable case is
given by Borisov and Fedorov \cite{BF}.
Let $R_2$, $\mathbf M$, $\mathbb I=\diag(I_1,I_2,I_3)$, be the radius, mass and the inertia operator of the ball $\mathrm B$,
 and let $R_1$ be the radius of the fixed sphere $\mathrm S$. There are three possible configurations:

\medskip
(i) rolling of $\mathrm B$ over outer surface of $\mathrm S$;

(ii) rolling of $\mathrm B$ over inner surface of $\mathrm S$ ($R_1>R_2$);

(iii) rolling of $\mathrm B$ over outer surface of $\mathrm S$, but $\mathrm S$ is within $\mathrm B$
($R_1<R_2$, in this case, the rolling ball $\mathrm B$ is
actually a spherical shell).

\medskip

Let
\begin{equation}\label{Epsilon}
\epsilon=\frac{R_1}{R_1\pm R_2},
\end{equation}
where we take "$+$" for the case (i) and "$-$" in the cases
(ii) and (iii) and let $D=\mathbf M R_2^2$.
The equations of motion in the frame attached to the ball can be written in the form
\begin{equation}
\dot{\mathbf G}={\mathbf G}\times\omega,
\qquad \dot{\gamma}=\epsilon
\gamma\times\omega, \label{Chap}\end{equation}
where
${\mathbf G}=\mathbb I \omega+ D\omega-D (\omega,\gamma)\gamma=\mathbf I\omega-D (\omega,\gamma)\gamma$  is the angular momentum of the ball with respect to the point of contact, and $\gamma$ is the unit
normal to the sphere $\mathrm S$ at the contact point. Here $\mathbf I=\mathbb I+D\mathbb E$, $\mathbb E=\diag(1,1,1)$.
When $R_1$ tends to infinity, $\epsilon$ tends to 1, $\gamma$ tends to the unit vector that is constant in the fixed reference frame.
This way we obtain the equations of motion of the Chaplygin ball rolling over the plane orthogonal to $\gamma$.

\medskip

In the space $\mathbb R^6({\omega},\gamma)$ the system has an invariant measure
with the density
\begin{equation}
\mu(\gamma)=\sqrt{(\gamma-D\mathbf I^{-1}\gamma,\gamma)}, \label{mu-ch}
\end{equation} the expression given by Chaplygin for $\epsilon=1$ \cite{Ch1}, and by Yaroshchuk for $\epsilon \ne 1$ \cite{Ya}. Also, the system \eqref{Chap} always has three integrals
\begin{equation}
F_1=(\gamma,\gamma)=1, \quad F_2=\frac12({\mathbf
G},\omega), \quad F_3=({\mathbf G},{\mathbf G}).
 \label{cl:int}
\end{equation}

For $\epsilon=1$, there is a fourth first integral $F_4=({\mathbf
G},\gamma)$. The problem is integrable by the Euler-Jacobi
theorem (see \cite{Kozlov2, AKN}): the phase space is almost everywhere foliated by
two-dimensional invariant tori with quasi-periodic, non-uniform motion
\eqref{Jacobi2} (see Chaplygin \cite{Ch1}). Moreover,
Borisov and Mamaev proved that the system \eqref{Chap} is
Hamiltonizable with respect to a certain nonlinear Poisson bracket
on $\mathbb R^6$ (\cite{BM}, see also
\cite{BM2008, BMT, BBM2013}).

\medskip

Remarkably, for $\epsilon=-1$ (the case (iii) with $R_2=2R_1$)
Borisov and Fedorov (see \cite{BF}) found an
integrable case with the following fourth first integral
\[
\tilde F_4=(I_2+I_3-I_1+D)\mathbf
G_1\gamma_1+(I_3+I_1-I_2+D)\mathbf
G_2\gamma_2+(I_1+I_2-I_3+D)\mathbf G_3\gamma_3.
\]
The system is
integrated on an invariant hypersurface $\tilde F_4^{-1}(0)$ \cite{BFM}. Its topological analysis
is given in \cite{BM3}.

\medskip

One can consider the additional nonholonomic constraint $(\omega,\gamma)=0$ describing \emph{no-twisting} condition: the ball $\mathrm B$ does not
rotate around the normal at the contact point (so called \emph{rubber Chaplygin ball}). Then the momentum with respect to the contact point can be expressed
as  $\mathbf G=\mathbb I\omega+D\omega=\mathbf I\omega$,
and the equations takes the form
\begin{equation}\label{e3a}
\dot{\mathbf G}=\mathbf G\times \omega+\lambda \gamma, \qquad
\dot{\gamma}=\epsilon \gamma\times\omega, \qquad (\omega,\gamma)=0,
\end{equation}
where the Lagrange multiplier
$
\lambda=-(\mathbf G,\mathbf I^{-1}(\mathbf G\times \omega))/(\gamma,\mathbf
I^{-1}\gamma)
$
is determined by differentiation of the constraint $(\omega,\gamma)=0$.
The system has an invariant measure with the density
$
\mu_\epsilon(\gamma)=(\mathbf I^{-1}\gamma,\gamma)^\frac{1}{2\epsilon}
$
(see \cite{EKMR2005} for $\epsilon=1$ and \cite{EK2007} for $\epsilon \ne 1$).
Apart of the integrability of the rolling over a
horizontal plane ($\epsilon=1$) \cite{EKMR2005}, as in the case of
non-rubber rolling, Borisov and Mamaev  proved the
integrability for $\epsilon=-1$ \cite{BM2007}. Note that for $\epsilon=1$, the above equations coincide with the equations of nonholonomic rigid body motion studied by Veselov and Veselova \cite{VeVe2}.

\medskip

The problem is Hamiltonizable for all $\epsilon$
\cite{EKMR2005, EK2007}. On the other hand, the rubber rolling of the ball where the mass
center does not coincide with the geometrical center over a
horizontal plane provides an example of the system having the
following interesting property (see \cite{BBM2}). The appropriate
phase space is foliated on invariant tori, such that the foliation
is isomorphic to the foliation of the integrable Euler case of the
rigid body motion about a fixed point, but the system itself has
not analytic invariant measure and is not Hamiltonizable.

\medskip

The gyroscopic generalizations of the mentioned Chaplygin ball problems are also well studied.
Markeev proved that the addition of a gyroscope to the Chaplygin ball problem of rolling of dynamically non-symmetric ball without slipping over a plane remains
integrable \cite{Markeev1985}.
As in Demchenko's thesis described in Sections \ref{bodyGyr} and \ref{ballGyr},
the addition of a gyroscope is equivalent to the addition of a constant angular momentum $\kappa$,
directed as the axis of the gyroscope, to $\mathbf G$ (with the new inertia operator described in Sections \ref{bodyGyr} and \ref{ballGyr}).
In the above notations, we can write the equations of the Chaplygin ball with the gyroscope rolling without slipping over the
\begin{equation}
\dot{\mathbf G}=({\mathbf G}+\kappa)\times\omega,
\qquad \dot{\gamma}=\epsilon
\gamma\times\omega. \label{Chap*}\end{equation}
When $\epsilon=1$ we have the Markeev integrable case \cite{Markeev1985}. The system has an invariant measure with the same density \eqref{mu-ch} and four first integrals
\begin{equation}
F_1=(\gamma,\gamma)=1, \quad F_2=\frac12({\mathbf
G},\omega), \quad F_3=({\mathbf G}+\kappa,{\mathbf G}+\kappa), \quad F_4=({\mathbf
G}+\kappa,\gamma)
 \label{cl:int*}
\end{equation}
The analysis of the bifurcation diagram and  the topology of the phase space of the Chaplygin ball with the gyroscope
case is studied in \cite{Moskvin2009} and \cite{Zhila2020}, respectively.

\medskip

The functions $F_1$, $F_2$, and $F_3$ are integrals for all $\epsilon$.
When the ball is dynamically symmetric with the gyroscope directed along the axis of the symmetry, it is the Bobylev--Zhukovsky case for $\epsilon=1$  \cite{Bob1892, Zhuk1893},
while when $\epsilon\ne 1$ and the Zhukovsky condition \eqref{uslovZh} on the moments of inertia of the ball and the gyroscope are satisfied, we obtain the Demchenko
integrable case. The integrability without the Zhukovskiy condition for the dynamically symmetric ball can be found  in  Borisov and Mamaev \cite{BoMa}.
Existence of an integrable case for a dynamically nonsymmetric ball with a gyroscope rolling over a sphere is still an open problem.

\medskip

The Voronec approach to the problem of rolling bodies given in \cite{Vor1911, Vor1912} can be found also in the recent papers \cite{Bychkov, KK, KZ}. In \cite{KK, KZ}, the problem of rolling without sliding of rotationally
symmetric  body on a fixed sphere is studied. It is assumed that the resultant of active forces is directed from the center of masses of the body to the center of a sphere. The problem reduces to
a linear differential equation of second order.  In a special case of motion of the nonhomogeneous dynamically symmetric ball, they proved the existence of Liouvillian solutions. In \cite{Bychkov}
a problem of rolling without slipping of a body with a gyroscope on a moving sphere is considered. It is assumed that the central ellipsoid of the system body + gyroscope is an ellipsoid of revolution. In a special case
when body is a sphere, the motion of the contact point is determined by quadratures. Analysis of trajectories of the contact point is given. This analysis, including pictures,
given in \cite{Bychkov} is very similar to the analysis presented by Demchenko in \cite{Dem1924}, see Figures \ref{fig:slika4}-\ref{fig:slika12}.

\medskip

Another line of current research is the application of
the Voronec equations \eqref{v4*}, \eqref{v5*} and their reductions in the case of symmetries \eqref{ChaplyginRed} to the study of their multi-dimensional versions, describing motions of the $n$--dimensional  ball rolling without slipping (and twisting) over a hyperplane or a sphere in $\mathbb R^n$ (see \cite{FK1995, Jo4, FNS, Jov2018b,  FNM, GajJov2019, Naranjo2019a, Naranjo2019b, Jov2019}). These examples, together with the classical one form a rich pool of nonholonomic systems. They motivate further study of the geometry and dynamics of nonholonomic systems including their integrabilty and Hamiltonization.

\medskip

There are very recent papers which are build on the results of Bilimovi\'c, e.g. \cite{BT2020, BMT2020}.
We hope that the current paper will further attract attention to the heritage of the Bilimovi\' c scientific school and their contribution to nonholonomic mechanical problems.
The Demchenko's integrable case and his comprehensive analysis provided in his doctoral thesis seem to be completely forgotten nowadays although still very modern and deserve to be known to a wider community.

\medskip

{\bf Acknowledgements.}
The authors would like to congratulate professor Veljko Vuji\v ci\'c, a distinguish member of the Bilimovi\'c school and the founding Editor-in-Chief of the Theoretical and Applied Mechanics, with his 90the anniversary.
This research has been partially supported by Mathematical Institute of the Serbian Academy of Sciences and Arts, the Science Fund of Serbia and the Ministry for Education, Science, and Technological Development of Serbia.

\end{document}